\documentclass[12pt,reqno]{amsart}

\usepackage[margin=1in]{geometry}
\usepackage{amsmath,amssymb,amsthm,mathtools}
\usepackage[dvipsnames]{xcolor}
\usepackage{booktabs,array}
\usepackage{tikz}
\usepackage{pgfplots}
\usepackage{float}
\pgfplotsset{compat=1.16}
\usetikzlibrary{decorations.pathreplacing,arrows.meta}
\usepackage[colorlinks=true,linkcolor=blue!60!black,citecolor=green!40!black,urlcolor=blue!60!black]{hyperref}

\newtheorem{theorem}{Theorem}[section]
\newtheorem{proposition}[theorem]{Proposition}
\newtheorem{lemma}[theorem]{Lemma}
\newtheorem{corollary}[theorem]{Corollary}
\theoremstyle{definition}
\newtheorem{definition}[theorem]{Definition}
\newtheorem{hypothesis}[theorem]{Hypothesis}
\newtheorem{conjecture}[theorem]{Conjecture}
\theoremstyle{remark}
\newtheorem{remark}[theorem]{Remark}

\numberwithin{equation}{section}

\newcommand{\R}{\mathbb{R}}
\newcommand{\N}{\mathcal{N}}
\newcommand{\G}{\mathcal{G}}
\newcommand{\ap}{\alpha_+}
\newcommand{\us}{u_\lambda^{*}}
\newcommand{\dd}{\,\mathrm{d}}
\newcommand{\yes}{$\checkmark$}
\newcommand{\no}{$\times$}

\begin{document}

\title[Bound states at Trudinger--Moser growth]{Bound states at Trudinger--Moser growth: obstructions, pairing and bubbling}

\author{Mohamed Majdoub}
\address{Department of Mathematics, College of Science, Imam Abdulrahman Bin Faisal University, P.O. Box 1982, Dammam 31441, Saudi Arabia}
\address{Basic and Applied Scientific Research Center, Imam Abdulrahman Bin Faisal University, P.O. Box 1982, Dammam 31441, Saudi Arabia}
\email{mmajdoub@iau.edu.sa} \email{med.majdoub@gmail.com}
\email{mohamed.majdoub@fst.rnu.tn}
\urladdr{https://orcid.org/0000-0001-6038-1069}

\begin{abstract}
We study radial solutions of the planar scalar field equation
$$
-\Delta u+u=\lambda u\,e^{u^2}, \qquad 0<\lambda<1,
$$
whose nonlinearity has critical Trudinger--Moser growth. These solutions describe standing waves of the two-dimensional nonlinear Schr\"odinger and Klein--Gordon equations. Motivated by Tang's recent uniqueness result for radial bound states with power nonlinearities in higher dimensions, we investigate the validity of the corresponding uniqueness theory at critical exponential growth. We show that the Serrin--Tang criterion and two structural hypotheses in Tang's analysis fail in this setting. We establish a pairing principle for the shooting problem, prove that Tang's variation conditions force the nodal count to be non-decreasing, and show that any two radial ground states intersect. We also derive a second-order expansion of the Liouville bubble arising at large initial heights and compute the corresponding flux correction. These asymptotics yield formal predictions for the first zero and the nodal structure at large heights. Our results suggest that the ground state is unique, whereas nodal radial bound states, if they exist, occur in pairs, contrary to the uniqueness phenomenon for power nonlinearities.
\end{abstract}

\maketitle

\section{Introduction}\label{sec:intro}

The nonlinear scalar field equation
\begin{equation}\label{eq:sfe}
\Delta u+g(u)=0\quad\text{in }\R^{n},\qquad u\in H^{1}(\R^{n}),\ u\not\equiv0,
\end{equation}
is the stationary problem behind the standing waves $\psi(t,x)=e^{i\omega t}u(x)$ of the nonlinear Schr\"odinger and Klein--Gordon equations. It has been one of the model problems of nonlinear analysis since the work of Strauss \cite{Strauss77} and of Berestycki and Lions \cite{BL83a,BL83b}. For $n\ge3$, Berestycki and Lions considered odd nonlinearities $g$ that are negative definite at the origin and at most Sobolev critical at infinity. They proved the existence of a positive radial ground state and of infinitely many radial solutions. The two-dimensional counterpart, in which the growth at infinity is subcritical in the exponential sense, is due to Berestycki, Gallou\"et and Kavian \cite{BGK83}. Radial solutions with any prescribed number of zeros were obtained by Jones and K\"upper \cite{JK86} and, by a direct argument for ordinary differential equations, by McLeod, Troy and Weissler \cite{MTW90}.

Uniqueness is much more delicate. For the model nonlinearity $g(u)=-u+|u|^{p-1}u$, the uniqueness of the positive solution was proved by Coffman \cite{Coffman72} for $n=3$ and $p=3$, then by McLeod and Serrin \cite{McLeodSerrin87} for a range of exponents, and finally by Kwong \cite{Kwong89} for all subcritical exponents in every dimension. We also refer to Peletier and Serrin \cite{PeletierSerrin83}, Chen and Lin \cite{ChenLin91}, Yanagida \cite{Yanagida91}, McLeod \cite{McLeod93}, Kabeya and Tanaka \cite{KabeyaTanaka99} and Serrin and Tang \cite{SerrinTang00}. The criterion of Serrin and Tang requires the function $u\,g'(u)/g(u)$ to be non-increasing beyond the positive zero of $g$. For uniqueness and nonuniqueness results in a more general radial setting, we refer to Ni and Nussbaum \cite{NiNussbaum85}. Uniqueness may indeed fail. For instance, D\'avila, del Pino and Guerra \cite{DdPG13} constructed three positive ground states in $\R^{3}$ for a suitable double-power nonlinearity. Positive solutions are radial by the theorem of Gidas, Ni and Nirenberg \cite{GNN81}; see also Li and Ni \cite{LiNi93}. The uniqueness of ground states is essential in dispersive theory. Indeed, it identifies the optimisers of the Gagliardo--Nirenberg inequality \cite{Weinstein83}, and it underlies the modulational stability theory \cite{Weinstein85}.

For sign-changing radial solutions, also called bound states or excited states, Berestycki and Lions conjectured that, at least for some classes of nonlinearities, there is exactly one radial solution with exactly $k$ nodes. Cort\'azar, Garc\'ia-Huidobro and Yarur \cite{CGHY11} proved the uniqueness of the bound state with one node for nonlinearities of the form $-|u|^{q-1}u+|u|^{p-1}u$ in low dimensions, including the plane. Ao, Wei and Yao \cite{AWY16} treated exponents close to the critical one. Cohen, Li and Schlag \cite{CLS24} proved by a computer-assisted argument the uniqueness of the first twenty excited states of $-\Delta u+u-u^{3}=0$ in $\R^{3}$, a question raised by Hastings and McLeod \cite{HastingsMcLeod12}. Multiplicity results for sign-changing bound states of other nonlinearities were obtained by Cort\'azar, Garc\'ia-Huidobro and Herreros \cite{CGHH15}. In a recent breakthrough, Tang \cite{Tang25} settled the conjecture for $g(u)=-u+|u|^{p-1}u$ with $1<p<(n+2)/(n-2)$, in every dimension $n\ge3$. His proof rests on a fine analysis of the variation $v=\partial u/\partial\alpha$ of the solution of the radial initial value problem with respect to the initial height $\alpha$. Tang's method has since been adapted to the logarithmic Schr\"odinger equation in dimensions $n\ge2$ \cite{LZS26}, and to sublinear equations, whose bound states have compact support \cite{ZZ26}. For the uniqueness of positive radial solutions of $m$-Laplacian equations in low dimensions, we refer to \cite{PZZ25}. The main uniqueness results quoted above are collected in Table~\ref{tab:known}.

\begin{table}[!htbp]
\centering
\caption{Uniqueness of radial solutions: a selection of known results. ``GS'' stands for ground state and ``$k$-BS'' for radial bound state with $k$ nodes. Entries marked with an asterisk are preprints.}\label{tab:known}
\small
\begin{tabular}{>{\raggedright\arraybackslash}p{4.4cm}>{\raggedright\arraybackslash}p{2.4cm}>{\raggedright\arraybackslash}p{3.3cm}>{\raggedright\arraybackslash}p{4.2cm}}
\toprule
Nonlinearity and setting & Dimension & Ground state & Bound states\\
\midrule
$-u+u^{3}$ in $\R^{3}$ & $n=3$ & Coffman \cite{Coffman72} & $k\le20$: Cohen--Li--Schlag \cite{CLS24}\\
$-u+|u|^{p-1}u$, subcritical & $n\ge3$ & Kwong \cite{Kwong89}; see also \cite{McLeodSerrin87,McLeod93} & all $k$: Tang \cite{Tang25}\\
$-u+|u|^{p-1}u$, $p>1$ & $n=2$ & Kwong \cite{Kwong89} & $k=1$: \cite{CGHY11}; all $k$: \cite{ZZ26}$^{*}$\\
$-|u|^{q-1}u+|u|^{p-1}u$ & $n\ge3$ (GS); $2\le n\le4$ ($k=1$) & Serrin--Tang \cite{SerrinTang00} & $k=1$: \cite{CGHY11}\\
$u\log u^{2}$ (logarithmic) & $n\ge2$ & see \cite{LZS26}$^{*}$ & all $k$: \cite{LZS26}$^{*}$\\
$u-|u|^{q-2}u$, $1<q<2$ (sublinear) & $n\ge2$ & see \cite{ZZ26}$^{*}$ & all $k$: \cite{ZZ26}$^{*}$\\
Critical exponential growth, Dirichlet problem in a disc & $n=2$ & \cite{Adimurthi98,Tang00,Tarsi03} & no radial nodal solution in small discs \cite{AY92}\\
$-u+\lambda ue^{u^{2}}$ in $\R^{2}$, problem \eqref{eq:P} & $n=2$ & existence \cite{AMS12,RufSani13}; uniqueness open (Conjecture~\ref{conj:ground}) & this paper: pairing (Theorem~\ref{thm:pairing}); conjectured non-uniqueness (Conjecture~\ref{conj:BL})\\
\bottomrule
\end{tabular}
\end{table}

In the plane every power $p>1$ is subcritical. Tang conjectures in the introduction of \cite{Tang25} that his theorem extends to $n=2$ and every $p>1$, a case to which, as he points out, his proof does not apply as written.

\begin{remark}\rm\label{rem:power}
The dimensional restriction enters the proof of \cite{Tang25} through the positivity of $Q_{n}=Q+nr^{n-1}u'v$ on the intervals $[b_{i},c_{i}]$, which is obtained in \cite[Lemma~5.3]{Tang25} from the Pohozaev bound $\omega(b_{i})>n-2$ when $n=3$, $p\ge2$ or $n\ge4$, and in \cite[Lemma~5.4]{Tang25} by an additional integral estimate when $n=3$ and $1<p<2$. When $n=2$, the function $Q_{n}$ coincides with $Q_{2}=Q+2ru'v$, and its positivity on the whole phase $[c_{i-1},c_{i}]$ follows directly from the induction hypothesis \cite[(5.8)]{Tang25}. Indeed, $Q_{2}'=2ru'v'>0$ on $(c_{i-1},\tau_{i})$ by \cite[(5.6)]{Tang25}, with $Q_{2}(c_{i-1})=Q(c_{i-1})>0$ (or $Q_{2}(0)=0$ when $i=1$), while $u'v>0$ on $(\tau_{i},c_{i})$ gives $Q_{2}>Q_{1}>0$ there, and $Q_{2}=Q_{1}>0$ at $\tau_{i}$ and at $c_{i}$. This observation, together with the modifications needed in the plane, is the subject of \cite[Section~6]{ZZ26}, where the uniqueness of the $k$-node bound states for $n=2$ and every $p>1$ is obtained as a by-product of the method of \cite{Tang25}. We do not use this result.
\end{remark}

The genuinely two-dimensional question, however, concerns a different notion of criticality. In $\R^{2}$, the space $H^{1}$ embeds into every $L^{q}$ with $2\le q<\infty$ but not into $L^{\infty}$, and the maximal growth compatible with a variational treatment is exponential. This is quantified by the Trudinger--Moser inequality \cite{Trudinger67,Moser71}. The supremum of $\int_{\Omega}e^{\beta u^{2}}$ over all $u\in H^{1}_{0}(\Omega)$ with $\|\nabla u\|_{2}\le1$ is finite if and only if $\beta\le4\pi$, and Carleson and Chang \cite{CarlesonChang86} proved that it is attained when $\Omega$ is a disc. Versions on the whole plane are due to Cao \cite{Cao92}, to Adachi and Tanaka \cite{AdachiTanaka00} and to Ruf \cite{Ruf05}. We also refer to Ishiwata \cite{Ishiwata11} for the existence and nonexistence of maximisers, and to Ibrahim, Masmoudi and Nakanishi \cite{IMN15} and Masmoudi and Sani \cite{MasmoudiSani15} for the inequality with the exact growth condition. A nonlinearity is said to have critical exponential growth if it behaves like $e^{4\pi s^{2}}$ at infinity \cite{dFMR95}. Elliptic problems in this regime were studied by Adimurthi \cite{Adimurthi90}, by Adimurthi and Yadava \cite{AY90}, by Cao \cite{Cao92} and by de Figueiredo, Miyagaki and Ruf \cite{dFMR95}. The corresponding loss of compactness was analysed by Adimurthi and Struwe \cite{AdimurthiStruwe00} and by Druet \cite{Druet06}, and, through profile decompositions of the critical embedding, by Bahouri, Majdoub and Masmoudi \cite{BMM11,BMM14}. For the scalar field equation on the whole plane with critical exponential growth, the existence of a ground state with a mountain-pass characterisation was proved by Alves, Montenegro and Souto \cite{AMS12} and by Ruf and Sani \cite{RufSani13,SaniThesis}, who followed the approach of Jeanjean and Tanaka \cite{JeanjeanTanaka03}.

Exponential nonlinearities are equally natural in dispersive equations. The Schr\"odinger equation in the Sobolev space of critical order was studied by Nakamura and Ozawa \cite{NakamuraOzawa98}. Ibrahim, Majdoub and Masmoudi \cite{IMM06} constructed global solutions of the two-dimensional Klein--Gordon equation with a nonlinearity of type $u(e^{4\pi u^{2}}-1)$ below the energy threshold dictated by the Trudinger--Moser inequality; see also the sharp double logarithmic inequality of \cite{IMM07}, the analysis of the energy supercritical regime in \cite{IMM11} and the uniqueness results of \cite{MajdoubMasmoudi15}. The energy-critical Schr\"odinger equation in the plane was treated by Colliander, Ibrahim, Majdoub and Masmoudi \cite{CIMM09}. Scattering for the energy-critical wave and Schr\"odinger equations with exponential nonlinearity was proved by Ibrahim, Majdoub, Masmoudi and Nakanishi \cite{IMMN09,IMMN12}. For the focusing Schr\"odinger equation with exponential nonlinearity, Dinh, Keraani and Majdoub \cite{DKM20} obtained thresholds for global existence and blow-up by variational arguments. We refer to \cite{BDM18,BDD20,BDM19} for weighted and time-oscillating exponential nonlinearities, and to \cite{IJM14,MOT18,MT21} for heat and biharmonic heat equations with exponential growth. The profiles of the standing waves of these equations are precisely the solutions of the planar scalar field equation with critical exponential growth (see Lemma~\ref{lem:dictionary} below). Their uniqueness and nondegeneracy are therefore the natural prerequisites for any modulation analysis near the ground state.

Uniqueness results at exponential growth are scarce, and the available ones concern the Dirichlet problem in a disc with a positive nonlinearity. Adimurthi \cite{Adimurthi98} proved the uniqueness of the positive solution of $-\Delta_{n} u=\lambda u^{n-1}e^{u}$ in the unit ball, which reduces to $-\Delta u=\lambda ue^{u}$ in the plane. His proof relies on the Emden--Fowler transformation used by Atkinson and Peletier \cite{AtkinsonPeletier86,AtkinsonPeletier86b}. Tang \cite{Tang00} obtained uniqueness for a large class of $n$-Laplacian Dirichlet problems in balls, which includes nonlinearities with exponential growth. Tarsi \cite{Tarsi03} combined the Emden--Fowler inversion with methods from the theory of sublinear equations and reached the exponents $1<\vartheta\le2$ in the nonlinearity $e^{u^{\vartheta}}$, including the critical case $\vartheta=2$. For sign-changing solutions, Adimurthi and Yadava \cite{AY92} proved that the Dirichlet problem in a disc with a critical exponential nonlinearity has no radial nodal solution when the radius is small. This is a two-dimensional analogue of the nonexistence result of Atkinson, Brezis and Peletier \cite{ABP90} for the Brezis--Nirenberg problem \cite{BrezisNirenberg83} in low dimensions. The blow-up of nodal radial solutions was analysed by Grossi and Naimen \cite{GrossiNaimen20}, and the nonexistence of least-energy nodal solutions for autonomous critical problems on $\R^{2}$ was shown by Alves and Pereira \cite{AlvesPereira15}.

\subsection{The problem and the results}
After the normalisation of Lemma~\ref{lem:dictionary}, we consider
\begin{equation}\label{eq:P}
-\Delta u+u=\lambda u\,e^{u^{2}}\quad\text{in }\R^{2},\qquad u\in H^{1}(\R^{2}),\qquad 0<\lambda<1. \tag{$P_\lambda$}
\end{equation}
We ask the two questions raised by Tang's theorem in this setting. Is the ground state of \eqref{eq:P} unique and nondegenerate? For each $k\ge1$, does \eqref{eq:P} have exactly one radial solution with exactly $k$ nodes?

The contribution of this paper is fourfold, and it is summarised in Table~\ref{tab:results}. First, we show in Section~\ref{sec:structure} that the Serrin--Tang criterion and two hypotheses of Tang's analysis fail at exponential growth for a structural reason (Lemma~\ref{lem:obstructions}). Second, we set up in Section~\ref{sec:shooting} the shooting framework for \eqref{eq:P} and prove a pairing principle (Theorem~\ref{thm:pairing}). Assume that the number $\N(\alpha)$ of zeros of the radial solution with initial height $\alpha$ exceeds $k$ for some $\alpha$ but not for large $\alpha$. Then \eqref{eq:P} has at least two distinct radial solutions, positive at the origin, with at most $k$ nodes. Under a natural simple-crossing property, they have exactly $k$ nodes, and under a transversality property the number of $k$-node bound states is either even or infinite (Proposition~\ref{prop:parity}). We also prove that Tang's two conditions on the variation imply that $\N$ is non-decreasing (Proposition~\ref{prop:reduction}). Consequently, they must fail near the right endpoint of every interval produced by the pairing principle (Corollary~\ref{cor:failure}). Moreover, any two radial ground states intersect (Lemma~\ref{lem:intersection}), which is the first step of a uniqueness argument by intersection. Third, we analyse in Section~\ref{sec:bubble} the Liouville bubble that forms at the origin as $\alpha\to\infty$. We prove the second-order expansion of the bubble in the core (Lemma~\ref{lem:core}) and compute exactly the correction of its flux (Lemma~\ref{lem:liouville}). We deduce from these facts, at a formal level, the asymptotic location \eqref{eq:z1} of the first zero and a mechanism which forces $\N(\alpha)=1$ for large $\alpha$. Fourth, we formulate in Section~\ref{sec:conjectures} precise conjectures, and we outline in Section~\ref{sec:program} a programme towards their proofs. Combined with the pairing principle, these conjectures say that $\N$ is \emph{not} monotone when $\lambda$ is close to $1$, that $k$-node bound states then occur in nested pairs, and that no nodal radial solution exists when $\lambda$ is below a threshold. In the language of Schr\"odinger standing waves, nodal radial standing waves would exist only at low frequency, and they would never be unique. Conjectures~\ref{conj:threshold} and~\ref{conj:BL} are also supported by numerical experiments, which are not reported here and are not used anywhere in the paper.

\begin{table}[H]
\centering
\caption{Main statements of the paper and their status.}\label{tab:results}
\small
\begin{tabular}{>{\raggedright\arraybackslash}p{8.3cm}>{\raggedright\arraybackslash}p{2.6cm}>{\raggedright\arraybackslash}p{3.6cm}}
\toprule
Statement & Status & Location\\
\midrule
Failure of the Serrin--Tang criterion and of the monotonicity of $g_{1}$ & proved & Lemma~\ref{lem:obstructions}\\
Semicontinuity and local constancy of the zero count $\N$ & proved & Lemmas~\ref{lem:local}, \ref{lem:small}\\
Pairing principle & proved & Theorem~\ref{thm:pairing}, Corollary~\ref{cor:pairs}\\
Exact number of nodes of the paired solutions & proved under Hypothesis~\ref{hyp:S} & Corollary~\ref{cor:exact}\\
Even number of $k$-node bound states & proved under transversality & Proposition~\ref{prop:parity}\\
Any two radial ground states intersect & proved & Lemma~\ref{lem:intersection}\\
Tang's conditions imply monotonicity and jumps of $\N$, hence fail near the paired heights & proved & Proposition~\ref{prop:reduction}, Corollary~\ref{cor:failure}\\
Core expansion of the bubble and flux correction $-4$ & proved & Lemmas~\ref{lem:liouville}, \ref{lem:core}\\
Asymptotics \eqref{eq:z1} of the first zero & formal & Section~\ref{sec:bubble}\\
$\N(\alpha)=1$ for large $\alpha$ & formal; Conjecture~\ref{conj:large} & Section~\ref{sec:bubble}\\
Uniqueness and nondegeneracy of the ground state & conjecture & Conjecture~\ref{conj:ground}\\
Threshold $\lambda_{c}$ for nodal solutions; pairs of $k$-node states & conjectures & Conjectures~\ref{conj:threshold}, \ref{conj:BL}\\
\bottomrule
\end{tabular}
\end{table}

\section{Setting and preliminaries}\label{sec:setting}

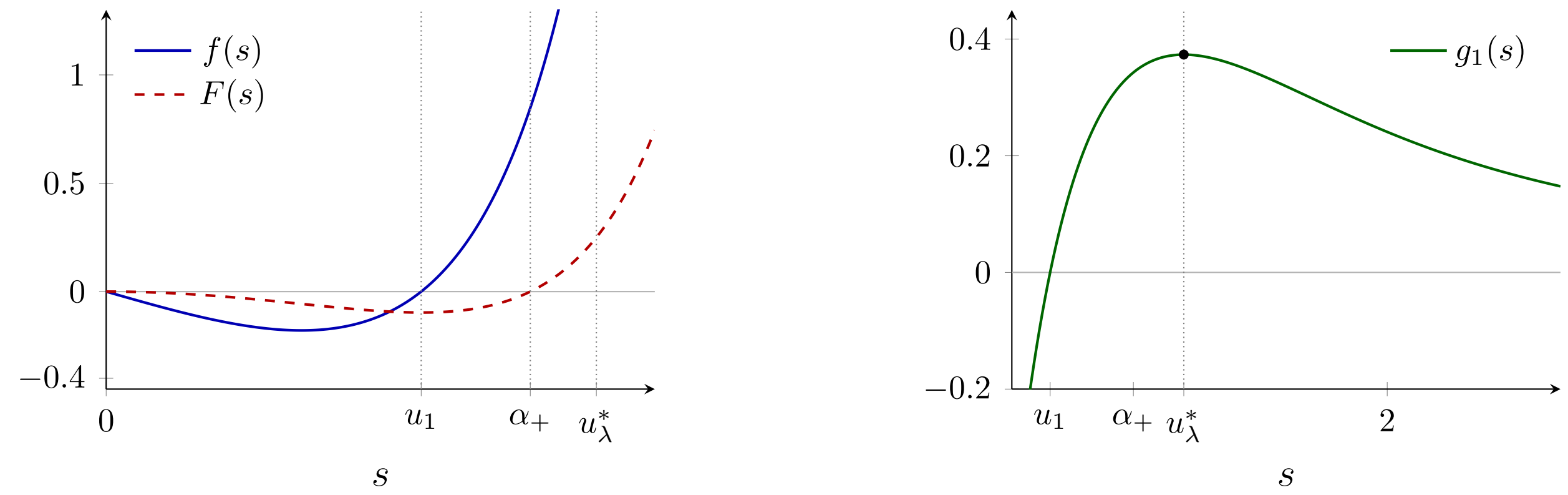
\begin{figure}[!htbp]
\centering
\begin{tikzpicture}
\begin{axis}[width=7.4cm,height=5.6cm,xmin=0,xmax=1.45,ymin=-0.45,ymax=1.3,axis lines=left,xlabel={$s$},samples=200,domain=0:1.45,legend style={at={(0.03,0.97)},anchor=north west,font=\footnotesize,draw=none,fill=none},xtick={0,0.8326,1.1209,1.2955},xticklabels={$0$,$u_{1}$,$\ap$,$\us$},ytick={-0.4,0,0.5,1},tick label style={font=\footnotesize}]
\draw[gray!60] (axis cs:0,0)--(axis cs:1.45,0);
\addplot[thick,blue!70!black] {-x+0.5*x*exp(x*x)};
\addlegendentry{$f(s)$}
\addplot[thick,dashed,red!70!black] {-x*x/2+0.25*(exp(x*x)-1)};
\addlegendentry{$F(s)$}
\draw[gray,densely dotted] (axis cs:0.8326,-0.45)--(axis cs:0.8326,1.3);
\draw[gray,densely dotted] (axis cs:1.1209,-0.45)--(axis cs:1.1209,1.3);
\draw[gray,densely dotted] (axis cs:1.2955,-0.45)--(axis cs:1.2955,1.3);
\end{axis}
\end{tikzpicture}\hfill
\begin{tikzpicture}
\begin{axis}[width=7.4cm,height=5.6cm,xmin=0.7,xmax=2.6,ymin=-0.2,ymax=0.45,axis lines=left,xlabel={$s$},samples=200,domain=0.7:2.6,legend style={at={(0.97,0.97)},anchor=north east,font=\footnotesize,draw=none,fill=none},xtick={0.8326,1.1209,1.2955,2},xticklabels={$u_{1}$,$\ap$,$\us$,$2$},ytick={-0.2,0,0.2,0.4},tick label style={font=\footnotesize}]
\draw[gray!60] (axis cs:0.7,0)--(axis cs:2.6,0);
\addplot[thick,green!40!black] {(0.5*exp(x*x)-1)/(0.5*x*x*exp(x*x))};
\addlegendentry{$g_{1}(s)$}
\draw[gray,densely dotted] (axis cs:1.2955,-0.2)--(axis cs:1.2955,0.45);
\fill (axis cs:1.2955,0.3734) circle (1.5pt);
\end{axis}
\end{tikzpicture}
\caption{The case $\lambda=\frac12$, for which $u_{1}\approx0.833$, $\ap\approx1.121$ and $\us\approx1.296$. Left: the nonlinearity $f$ and its primitive $F$, which vanish at $u_{1}$ and $\ap$ respectively. Right: the function $g_{1}$ of Lemma~\ref{lem:obstructions}(b), which increases up to $\us$ and then decreases to $0$, in contrast with the power case where $g_{1}$ is increasing on $(0,\infty)$.}\label{fig:levels}
\end{figure}

Throughout the paper, $0<\lambda<1$ and
\begin{equation}\label{eq:f}
f(s)=-s+\lambda s\,e^{s^{2}},\qquad F(s)=\int_{0}^{s}f(\sigma)\dd\sigma=-\frac{s^{2}}{2}+\frac{\lambda}{2}\big(e^{s^{2}}-1\big),
\end{equation}
so that \eqref{eq:P} reads $\Delta u+f(u)=0$. We write $\zeta=1-\lambda=-f'(0)$, and we introduce the three positive numbers
\begin{equation}\label{eq:levels}
u_{1}=\sqrt{\log(1/\lambda)},\qquad \ap=\sqrt{s_{+}},\qquad \us=\sqrt{s^{*}},
\end{equation}
where $s_{+}>0$ is the positive root of $\lambda(e^{s}-1)=s$ and $s^{*}>0$ is the root of $(1+s)e^{-s}=\lambda$. The existence and uniqueness of $s_{+}$ and $s^{*}$ are proved in Lemma~\ref{lem:elementary} and Lemma~\ref{lem:obstructions} respectively. These levels are illustrated in Figure~\ref{fig:levels}.

\begin{lemma}\label{lem:dictionary}
~\begin{enumerate}
\item[(i)] A function $u$ solves \eqref{eq:P} if and only if $w=u/\sqrt{4\pi}$ solves $-\Delta w+w=\lambda w\,e^{4\pi w^{2}}$ in $\R^{2}$, which is the equation studied in \cite{RufSani13,SaniThesis}.
\item[(ii)] Let $\omega>0$. The function $\psi(t,x)=e^{i\omega t}\varphi(x)$, with $\varphi$ real valued, solves the focusing Schr\"odinger equation $i\partial_{t}\psi+\Delta\psi+(e^{4\pi|\psi|^{2}}-1)\psi=0$ if and only if $u(x)=\sqrt{4\pi}\,\varphi\big(x/\sqrt{1+\omega}\big)$ solves \eqref{eq:P} with $\lambda=1/(1+\omega)$.
\item[(iii)] Let $|\omega|<1$. The function $\phi(t,x)=e^{i\omega t}\varphi(x)$, with $\varphi$ real valued, solves the focusing Klein--Gordon equation $\partial_{t}^{2}\phi-\Delta\phi+\phi=(e^{4\pi|\phi|^{2}}-1)\phi$ if and only if $u(x)=\sqrt{4\pi}\,\varphi\big(x/\sqrt{2-\omega^{2}}\big)$ solves \eqref{eq:P} with $\lambda=1/(2-\omega^{2})\in[\tfrac12,1)$.
\end{enumerate}
\end{lemma}

\begin{proof}
For (i), the substitution $u=\sqrt{4\pi}\,w$ gives $\lambda u e^{u^{2}}=\sqrt{4\pi}\,\lambda w e^{4\pi w^{2}}$, and the claim follows after division by $\sqrt{4\pi}$. For (ii), the time derivative of $e^{i\omega t}\varphi$ produces the term $-\omega\varphi$, and we obtain the standing-wave equation $-\Delta\varphi+(1+\omega)\varphi=\varphi e^{4\pi\varphi^{2}}$. We write $\varphi(y)=(4\pi)^{-1/2}u(\sqrt{1+\omega}\,y)$. Then $\Delta_{y}\varphi=(4\pi)^{-1/2}(1+\omega)(\Delta u)(\sqrt{1+\omega}\,y)$ and $4\pi\varphi^{2}=u^{2}$, so the equation becomes $(1+\omega)(-\Delta u+u)=ue^{u^{2}}$, which is \eqref{eq:P} with $\lambda=1/(1+\omega)$. The proof of (iii) is identical, since the standing-wave equation now reads $-\Delta\varphi+(2-\omega^{2})\varphi=\varphi e^{4\pi\varphi^{2}}$. Finally, the range of $\lambda$ follows from $0\le\omega^{2}<1$.
\end{proof}

\begin{remark}\rm\label{rem:mu1}
The ground states relevant for the Klein--Gordon equation correspond to $\lambda\in[\frac12,1)$, while the Schr\"odinger standing waves cover the whole range $\lambda\in(0,1)$, high frequencies corresponding to small values of $\lambda$. The Schr\"odinger nonlinearity in Lemma~\ref{lem:dictionary}(ii) is the case $\mu=0$ of the family $(e^{4\pi|u|^{2}}-1-4\pi\mu|u|^{2})u$, $\mu\in\{0,1\}$, studied in \cite{DKM20}. For $\mu=1$, the same computation leads to $\Delta u+f_{1}(u)=0$ with $f_{1}(s)=s\big({-1}+\lambda(e^{s^{2}}-1-s^{2})\big)$ and $\lambda=1/\omega\in(0,\infty)$. The function $f_{1}(s)/s$ is strictly increasing on $(0,\infty)$, $f_{1}'(0)=-1$, and $F_{1}(s)=\frac12\Phi_{1}(s^{2})$ with $\Phi_{1}(\sigma)=\lambda(e^{\sigma}-1-\sigma-\sigma^{2}/2)-\sigma$, which is strictly convex with $\Phi_{1}(0)=0$ and $\Phi_{1}'(0)=-1$. Hence $f_{1}$ satisfies the conditions (C1)--(C4) recalled in Section~\ref{sec:shooting}, and the rigorous results of Section~\ref{sec:shooting}, whose proofs use only these conditions and the monotonicity of $f(s)/s$, hold verbatim for $f_{1}$. We do not pursue this family further.
\end{remark}

\begin{lemma}\label{lem:elementary}
The function $f$ is smooth and odd, and the following properties hold.
\begin{enumerate}
\item[(a)] $f<0$ on $(0,u_{1})$, $f>0$ on $(u_{1},\infty)$, $f'(0)=-\zeta<0$ and $f'(u_{1})=2u_{1}^{2}>0$.
\item[(b)] $F<0$ on $(0,\ap)$ and $F>0$ on $(\ap,\infty)$, where $\ap$ is well defined by \eqref{eq:levels}. Moreover, $\ap>u_{1}$.
\item[(c)] $sf(s)-2F(s)=\lambda\big(s^{2}e^{s^{2}}-e^{s^{2}}+1\big)>0$ for $s\neq0$.
\item[(d)] $f''(s)=\lambda s\,e^{s^{2}}(6+4s^{2})>0$ and $sf'(s)-f(s)=2\lambda s^{3}e^{s^{2}}>0$ for $s>0$.
\item[(e)] The function $s\mapsto f(s)/s=\lambda e^{s^{2}}-1$ is strictly increasing on $(0,\infty)$.
\end{enumerate}
\end{lemma}

\begin{proof}
Since $f(s)=s(\lambda e^{s^{2}}-1)$, the sign of $f$ on $(0,\infty)$ is that of $\lambda e^{s^{2}}-1$, which changes exactly at $u_{1}$. Moreover, $f'(s)=-1+\lambda e^{s^{2}}(1+2s^{2})$, hence $f'(0)=\lambda-1$, and the identity $\lambda e^{u_{1}^{2}}=1$ yields $f'(u_{1})=2u_{1}^{2}$. This proves (a), and (e) follows immediately from the same factorisation. For (b), we observe that $F(s)=\frac12\Phi(s^{2})$ with $\Phi(\sigma)=\lambda(e^{\sigma}-1)-\sigma$. The function $\Phi$ is strictly convex, vanishes at $0$, satisfies $\Phi'(0)=\lambda-1<0$, and tends to $+\infty$ at infinity. We infer that $\Phi$ has exactly one positive zero $s_{+}$, and that it is negative on $(0,s_{+})$ and positive on $(s_{+},\infty)$. Since $\Phi'(\sigma)=\lambda e^{\sigma}-1<0$ on $(0,\log(1/\lambda))$, the function $\Phi$ is negative on $(0,\log(1/\lambda)]$, and we conclude that $s_{+}>\log(1/\lambda)$, that is, $\ap>u_{1}$. For (c), a direct computation gives $sf(s)-2F(s)=\lambda\chi(s^{2})$ with $\chi(\sigma)=\sigma e^{\sigma}-e^{\sigma}+1$. As $\chi(0)=0$ and $\chi'(\sigma)=\sigma e^{\sigma}>0$ for $\sigma>0$, the expression is positive for $s\neq0$. Finally, (d) follows by differentiating $f'$ once more, and from the identity $sf'(s)=f(s)+2\lambda s^{3}e^{s^{2}}$.
\end{proof}

\begin{remark}\rm\label{rem:2d}
In dimension two, the Pohozaev function of the radial equation degenerates. We use the notation of \cite{Tang25}, where the nonlinearity is normalised so that its positive zero equals $1$ (see the paragraph before Proposition~\ref{prop:tang}). The functions entering conditions (C5) and (C7) of \cite{Tang25} are $h(s)=2nF(s)-(n-2)sf(s)$ and $h_{2}(s)=(n+2)f(s)-(n-2)sf'(s)$. When $n=2$, we get $h=4F$ and $h_{2}=4f$. Consequently, condition (C5) holds with threshold equal to the positive zero of $F$ (in our normalisation, $\ap$), whereas for $n\ge3$ the threshold in (C5) is strictly larger than the positive zero of $F$. Likewise, condition (C7) holds with threshold equal to the positive zero of $f$ (in our normalisation, $u_{1}$). In addition, the Pohozaev function $P(r)=r^{2}\big(u'(r)^{2}+2F(u(r))\big)$ used in Lemma~\ref{lem:pohozaev} below reduces to $P(r)=2r^{2}E(r)$, where $E$ is the energy \eqref{eq:energy}.
\end{remark}

\section{Structural obstructions}\label{sec:structure}

Let $K(s)=sf'(s)/f(s)$ for $s>u_{1}$. The uniqueness criterion of Serrin and Tang \cite{SerrinTang00} requires $K$ to be non-increasing on $(u_{1},\infty)$. Condition (C8) of \cite{Tang25} requires $f/(sf')=1/K$ to be increasing on $(u_{1},\infty)$; when $f'>0$ there, this is equivalent to the Serrin--Tang condition. In the analysis of the first phase in \cite[Lemma~4.2]{Tang25}, the function
\[
g_{1}(s)=\frac{2f(s)}{sf'(s)-f(s)}
\]
is assumed to satisfy $g_{1}'>0$ on $(0,\infty)$. The functions $K$, $1/K$ and $g_{1}$ are invariant under the normalisation $w=u/u_{1}$ of Section~\ref{sec:shooting}, in the sense that the corresponding functions of $\tilde f(w)=f(u_{1}w)/u_{1}$ are obtained by evaluating them at $s=u_{1}w$. We may therefore state all conditions directly for $f$.

\begin{lemma}\label{lem:obstructions}
~\begin{enumerate}
\item[(a)] For every $s>u_{1}$ one has $K(s)>1+2s^{2}$. In particular, $K$ is unbounded, it is not non-increasing on any half-line $(b,\infty)$, and $1/K$ is not increasing on $(u_{1},\infty)$.
\item[(b)] For $s>0$,
\[
g_{1}(s)=\frac{\lambda e^{s^{2}}-1}{\lambda s^{2}e^{s^{2}}},\qquad g_{1}'(s)=\frac{2}{s^{3}}\left(\frac{(1+s^{2})e^{-s^{2}}}{\lambda}-1\right).
\]
Hence $g_{1}$ is strictly increasing on $(0,\us]$ and strictly decreasing on $[\us,\infty)$, where $\us$ is well defined by \eqref{eq:levels}. Moreover, $\ap<\us$.
\end{enumerate}
\end{lemma}

\begin{proof}
Let $s>u_{1}$. By Lemma~\ref{lem:elementary}(d), $sf'(s)=f(s)+2\lambda s^{3}e^{s^{2}}$, and $f(s)=s(\lambda e^{s^{2}}-1)>0$. We obtain
\[
K(s)=1+\frac{2s^{2}\,\lambda e^{s^{2}}}{\lambda e^{s^{2}}-1}>1+2s^{2},
\]
because $\lambda e^{s^{2}}>\lambda e^{s^{2}}-1>0$. A function which is non-increasing on $(b,\infty)$ is bounded above there by its value at any point of that interval, and this is incompatible with $K(s)\to\infty$ as $s\to\infty$. Likewise, the positive function $1/K$ tends to $0$ at infinity, and therefore it cannot be increasing on $(u_{1},\infty)$. This proves (a).

For (b), the formula for $g_{1}$ follows from Lemma~\ref{lem:elementary}(d). Writing $g_{1}(s)=s^{-2}-\lambda^{-1}s^{-2}e^{-s^{2}}$, we find
\[
g_{1}'(s)=-\frac{2}{s^{3}}+\frac{1}{\lambda}\Big(\frac{2}{s^{3}}+\frac{2}{s}\Big)e^{-s^{2}}=\frac{2}{s^{3}}\left(\frac{(1+s^{2})e^{-s^{2}}}{\lambda}-1\right).
\]
The function $\sigma\mapsto(1+\sigma)e^{-\sigma}$ has derivative $-\sigma e^{-\sigma}<0$ on $(0,\infty)$, and it decreases from $1$ to $0$. Since $0<\lambda<1$, it takes the value $\lambda$ exactly once, at $s^{*}$, and it is larger than $\lambda$ on $(0,s^{*})$ and smaller on $(s^{*},\infty)$. The sign of $g_{1}'$ is therefore as stated. Finally, the definition of $s_{+}$ gives $\lambda e^{s_{+}}=s_{+}+\lambda$, so that
\[
(1+s_{+})e^{-s_{+}}=\frac{\lambda(1+s_{+})}{s_{+}+\lambda}>\lambda,
\]
since $\lambda<1$. By the monotonicity just established, we deduce that $s_{+}<s^{*}$, that is, $\ap<\us$.
\end{proof}

\begin{remark}\rm
Lemma~\ref{lem:obstructions}(a) is the precise sense in which the Serrin--Tang approach, designed for nonlinearities of at most power growth, does not reach Trudinger--Moser growth. By part (b), the monotonicity of $g_{1}$ required in \cite[Lemma~4.2]{Tang25} fails on $(\us,\infty)$, and it holds on the range $(0,\alpha]$ of a positive solution only if $\alpha\le\us$. In the proof of \cite[Lemma~4.2]{Tang25}, the derivative $g_{1}'$ only enters through the identity $T_{1}'(r)=-g_{1}'(u(r))u'(r)M(r)$, and it is evaluated at the values $u(r)\in(0,\alpha)$ taken for $0<r<z_{1}$ (or for all $r>0$ in the case of a ground state). The remaining hypotheses of that lemma, namely (C1), (C3), (C5), (C6), (C9) and the positivity of $sf'(s)-f(s)$, hold for our $f$ in dimension two by Lemma~\ref{lem:elementary} and Remark~\ref{rem:2d}. Since \cite[Sections~2 and~4]{Tang25} are written for a real dimension $n>1$, we infer that the conclusions of \cite[Lemma~4.2]{Tang25} hold for every ground state or nodal solution of \eqref{eq:ivp} whose initial height satisfies $\alpha\le\us$. Every positive solution which decays at infinity has $\alpha>\ap$ (Proposition~\ref{prop:tang}(ii)), and $\ap<\us$, so the condition $\alpha\le\us$ is not excluded a priori. However, a uniqueness proof must handle all ground states simultaneously, and the condition $\alpha\le\us$ cannot be verified for an unknown second ground state.
\end{remark}

Table~\ref{tab:conditions} compares the structural conditions of \cite{Tang25} for the power nonlinearity and for the nonlinearity of \eqref{eq:P}.

\begin{table}[!htbp]
\centering
\caption{Structural conditions of \cite{Tang25} for the power nonlinearity $-u+|u|^{p-1}u$ in $\R^{n}$, $n\ge3$, and for the nonlinearity $f$ of \eqref{eq:P} in $\R^{2}$. All conditions are stated after normalising the positive zero of the nonlinearity.}\label{tab:conditions}
\small
\begin{tabular}{>{\raggedright\arraybackslash}p{5.6cm}cc>{\raggedright\arraybackslash}p{4.3cm}}
\toprule
Condition & Power, $n\ge3$ & $f$ of \eqref{eq:P}, $n=2$ & Justification for $f$\\
\midrule
(C1)--(C4): sign, primitive, $f'(0)<0$, $f'>0$ at the zero & \yes & \yes & Lemma~\ref{lem:elementary}(a)--(b)\\
(C5): sign of $h=2nF-(n-2)sf$ & \yes & \yes & $h=4F$, Remark~\ref{rem:2d}\\
(C6): $sf-2F\ge0$ & \yes & \yes & Lemma~\ref{lem:elementary}(c)\\
(C7): sign of $h_{2}=(n+2)f-(n-2)sf'$ & \yes & \yes & $h_{2}=4f$, Remark~\ref{rem:2d}\\
(C8): $f/(sf')$ increasing (Serrin--Tang) & \yes & \no & Lemma~\ref{lem:obstructions}(a)\\
(C9): $f''>0$ beyond the zero of $f$ & \yes & \yes & Lemma~\ref{lem:elementary}(d)\\
$sf'-f>0$ on $(0,\infty)$ & \yes & \yes & Lemma~\ref{lem:elementary}(d)\\
$g_{1}'>0$ on $(0,\infty)$ & \yes & \no\ (only on $(0,\us)$) & Lemma~\ref{lem:obstructions}(b)\\
Connection identities \cite[Appendix~A.4]{Tang25} & \yes & not available & specific to the power case\\
\bottomrule
\end{tabular}
\end{table}

\section{The shooting framework and rigorous results}\label{sec:shooting}

Radial solutions $u(x)=u(r)$, $r=|x|$, of \eqref{eq:P} are obtained from the initial value problem
\begin{equation}\label{eq:ivp}
u''+\frac{1}{r}u'+f(u)=0,\quad r>0,\qquad u(0)=\alpha>0,\quad u'(0)=0.
\end{equation}
Problem \eqref{eq:ivp} is equivalent to the Volterra integral equation
\[
u(r)=\alpha-\int_{0}^{r}\tau\log\Big(\frac{r}{\tau}\Big)f(u(\tau))\dd\tau,
\]
whose kernel is bounded on $\{0\le\tau\le r\le R\}$. Since $f$ is smooth, a contraction argument gives a unique local solution $u(\cdot,\alpha)$. Moreover, the solution depends smoothly on $\alpha$ in $C^{1}([0,R])$ for every $R$ for which it exists \cite{BLP81,Hartman02,Tang25}. The energy
\begin{equation}\label{eq:energy}
E(r)=\tfrac12u'(r)^{2}+F(u(r)),\qquad E'(r)=-\frac{u'(r)^{2}}{r},
\end{equation}
is non-increasing. As $F$ is even and coercive, the inequality $F(u(r))\le E(r)\le F(\alpha)$ bounds $u$ on the interval of existence. Then $\frac12u'(r)^{2}\le F(\alpha)-\min F$ bounds $u'$, and we conclude that the solution is global. The variation $v=\partial u/\partial\alpha$ is well defined for all $r\ge0$, and it solves
\[
v''+\frac{1}{r}v'+f'(u)v=0,\qquad v(0)=1,\quad v'(0)=0.
\]

\begin{remark}\rm\label{rem:regularity}
Every solution $u\in H^{1}(\R^{2})$ of \eqref{eq:P} is smooth, bounded, and tends to zero at infinity. Indeed, $\int_{\R^{2}}(e^{\beta u^{2}}-1)<\infty$ for every $\beta>0$ by \cite{Cao92,Ruf05}. Hence, for every $q<\infty$, the integral $\int_{B_{1}(x)}e^{qu^{2}}$ is bounded independently of $x\in\R^{2}$. Since $|s|\le e^{s^{2}}$, the function $-u+\lambda ue^{u^{2}}$ is bounded in $L^{q}(B_{1}(x))$ independently of $x$. The interior $L^{q}$ estimates then show that $u$ is bounded in $W^{2,q}(B_{1/2}(x))$ independently of $x$. We infer that $u$ is bounded and uniformly continuous on $\R^{2}$, and a bootstrap argument shows that $u$ is smooth. As $u\in L^{2}(\R^{2})$ is uniformly continuous, it tends to zero at infinity. In particular, every radial solution in $H^{1}(\R^{2})$ with $u(0)=\alpha>0$ coincides with $u(\cdot,\alpha)$. Moreover, every positive solution is radially symmetric about some point by \cite{GNN81,LiNi93}, since $f(s)=-\zeta s+O(s^{3})$ as $s\to0$. These facts are not used below.
\end{remark}

\begin{definition}\label{def:types}
The solution $u=u(\cdot,\alpha)$ of \eqref{eq:ivp} is called a \emph{ground state} if $u>0$ and $u(r)\to0$ as $r\to\infty$. It is called a \emph{$k$-node bound state} if it has exactly $k\ge1$ zeros in $(0,\infty)$ and $u(r)\to0$ as $r\to\infty$. It is called \emph{oscillatory} if $E(r)\le0$ for some $r\ge0$. We denote by $\N(\alpha)$ the number of zeros of $u(\cdot,\alpha)$ in $(0,\infty)$, and by $\G$ the set of all $\alpha>0$ for which $u(\cdot,\alpha)$ is a ground state or a bound state.
\end{definition}

The following proposition collects results of Tang \cite[Propositions~2.1--2.4]{Tang25}; see also \cite{BLP81,PeletierSerrin83}. They are proved in \cite{Tang25} for a general nonlinearity satisfying conditions (C1)--(C4) there, namely: $f\in C^{1}(\R)$ is odd, $f<0$ on $(0,1)$ and $f>0$ on $(1,\infty)$; $F$ has a positive zero $\alpha^{*}>1$; $f'(0)<0$; and $f'(1)>0$. In \cite[Section~2]{Tang25}, the dimension $n$ appearing in the radial equation is treated as a real parameter in $(1,\infty)$, so these results apply when $n=2$. Tang normalises the positive zero of $f$ to be $1$, and this is inessential. Indeed, $w=u/u_{1}$ solves $w''+w'/r+\tilde f(w)=0$ with $\tilde f(w)=f(u_{1}w)/u_{1}$. By Lemma~\ref{lem:elementary}(a)--(b), the function $\tilde f$ satisfies (C1)--(C4), with $\tilde f'(0)=-\zeta$, $\tilde f'(1)=2u_{1}^{2}$ and $\alpha^{*}=\ap/u_{1}$.

\begin{proposition}[\cite{Tang25}]\label{prop:tang}
Let $u=u(\cdot,\alpha)$ with $\alpha>0$ and $\alpha\neq u_{1}$ (for $\alpha=u_{1}$ we have $u\equiv u_{1}$).
\begin{enumerate}
\item[(i)] $E$ is strictly decreasing on $(0,\infty)$. If $\alpha>u_{1}$, then $\|u\|_{\infty}=\alpha$.
\item[(ii)] If $u$ is a ground state or has a zero, then $\alpha>\ap$. If $\alpha\le\ap$, then $u>0$ and $u$ is oscillatory.
\item[(iii)] Assume that $u$ has zeros $z_{1}<z_{2}<\cdots$. Then $u$ has exactly one critical point $c_{i}\in(z_{i},z_{i+1})$ between two consecutive zeros, $|u(c_{i})|>\ap$, and $\alpha>|u(c_{1})|>|u(c_{2})|>\cdots$. Moreover, $\N(\alpha)<\infty$.
\item[(iv)] Exactly one of the following occurs: $u$ is a ground state, $u$ is a bound state, or $u$ is oscillatory. Ground states and bound states satisfy $E>0$ on $(0,\infty)$ and $u'(r)/u(r)\to-\sqrt{\zeta}$ as $r\to\infty$. Moreover, $\limsup_{r\to\infty}|u(r)|e^{(\sqrt{\zeta}-\varepsilon)r}<\infty$ for every $\varepsilon\in(0,\sqrt\zeta)$, and the same exponential decay holds for $u'$. An oscillatory solution satisfies $E(r)<0$ for all large $r$, has no zero beyond the first $r$ at which $E(r)\le0$, and oscillates about $u_{1}$ or $-u_{1}$.
\end{enumerate}
\end{proposition}

Several parts of Proposition~\ref{prop:tang} are elementary consequences of \eqref{eq:energy}. For instance, at a zero $z$ of $u$ we have $E(z)=\frac12u'(z)^{2}>0$, the zero being simple by uniqueness for \eqref{eq:ivp}. Since $E$ is strictly decreasing, we get $F(\alpha)=E(0)>E(z)>0$, and Lemma~\ref{lem:elementary}(b) yields $\alpha>\ap$. In particular, every $\alpha\in\G$ yields a classical, exponentially decaying, and hence finite-energy solution of \eqref{eq:P}.

\begin{lemma}[Pohozaev identity in the plane]\label{lem:pohozaev}
If $\alpha\in\G$, then $\int_{0}^{\infty}rF(u(r,\alpha))\dd r=0$.
\end{lemma}

\begin{proof}
Let $P(r)=r^{2}\big(u'^{2}+2F(u)\big)=2r^{2}E(r)$. Using \eqref{eq:energy}, we get $P'(r)=4rE(r)-2ru'(r)^{2}=4rF(u(r))$. We have $P(0)=0$. Moreover, $P(r)\to0$ as $r\to\infty$ by the exponential decay of $u$ and $u'$, since $F(s)=O(s^{2})$ as $s\to0$. Integrating over $(0,\infty)$, we obtain the claim.
\end{proof}

\subsection{Semicontinuity and local constancy of the zero count}

\begin{lemma}\label{lem:local}
Let $\alpha>0$.
\begin{enumerate}
\item[(a)] There is $\delta>0$ such that $\N(\beta)\ge\N(\alpha)$ for $|\beta-\alpha|<\delta$.
\item[(b)] If $\alpha\notin\G$, there is $\delta>0$ such that $\N(\beta)=\N(\alpha)$ for $|\beta-\alpha|<\delta$.
\end{enumerate}
\end{lemma}

\begin{proof}
Every zero $z>0$ of $u(\cdot,\alpha)$ is simple. Indeed, $u(z)=u'(z)=0$ would force $u\equiv0$ on $(0,\infty)$ by uniqueness for \eqref{eq:ivp} away from the origin, and then $u(0)=0$ by continuity, which contradicts $u(0)=\alpha>0$. For (a), we may assume $\N(\alpha)\ge1$; recall that $\N(\alpha)$ is finite by Proposition~\ref{prop:tang}(iii). We choose $R$ larger than the last zero of $u(\cdot,\alpha)$, so that $u(R,\alpha)\neq0$. As $\beta\to\alpha$, we have $u(\cdot,\beta)\to u(\cdot,\alpha)$ in $C^{1}([0,R])$. Since $u(\cdot,\alpha)$ changes sign at each of its zeros, the intermediate value theorem gives a zero of $u(\cdot,\beta)$ near each zero of $u(\cdot,\alpha)$ in $(0,R)$. We conclude that $u(\cdot,\beta)$ has at least $\N(\alpha)$ zeros when $\beta$ is close to $\alpha$.

For (b), suppose first that $\alpha=u_{1}$. Then $E(0)=F(u_{1})<0$ by Lemma~\ref{lem:elementary}(b), so that $F(\beta)<0$ for $\beta$ close to $u_{1}$. Since the energy is non-increasing, the solution $u(\cdot,\beta)$ satisfies $E\le F(\beta)<0$ on $[0,\infty)$. As $E=\frac12u'^{2}\ge0$ at any zero, it follows that $u(\cdot,\beta)$ has no zero, and $\N\equiv0$ near $u_{1}$.

Suppose now that $\alpha\neq u_{1}$ and $\alpha\notin\G$. By Proposition~\ref{prop:tang}(iv), $u(\cdot,\alpha)$ is oscillatory, and the strict monotonicity of $E$ gives some $\bar r>0$ with $E(\bar r)<0$. Then $u(\bar r,\alpha)\neq0$, and all zeros of $u(\cdot,\alpha)$ lie in $(0,\bar r)$. Let $z_{1}<\dots<z_{m}$ be these zeros, with $m=\N(\alpha)$. Since $u'(z_{i},\alpha)\neq0$, we may choose disjoint closed intervals $I_{i}\subset(0,\bar r)$ centred at $z_{i}$ on which $|u'(\cdot,\alpha)|$ is bounded below by a positive constant. Then $|u(\cdot,\alpha)|$ is bounded below by a positive constant on the compact set $[0,\bar r]\setminus\bigcup_{i}\mathrm{int}(I_{i})$. By the $C^{1}$ convergence on $[0,\bar r]$, for $\beta$ close to $\alpha$ the function $u(\cdot,\beta)$ does not vanish on this set, it is strictly monotone on each $I_{i}$, and it changes sign on each $I_{i}$. We infer that $u(\cdot,\beta)$ has exactly one zero in each $I_{i}$, hence exactly $m$ zeros in $(0,\bar r]$. Moreover, its energy at $\bar r$ is still negative, by continuity. Since the energy is non-increasing, $u(\cdot,\beta)$ has no zero in $[\bar r,\infty)$, and we conclude that $\N(\beta)=\N(\alpha)$.
\end{proof}

\begin{lemma}\label{lem:small}
$\N(\beta)=0$ for every $\beta\in(0,\ap]$, and there is $\delta_{0}>0$ such that $\N=0$ on $(0,\ap+\delta_{0})$.
\end{lemma}

\begin{proof}
For $\beta\in(0,\ap]$ with $\beta\neq u_{1}$, Proposition~\ref{prop:tang}(ii) gives $u(\cdot,\beta)>0$, and for $\beta=u_{1}$ the solution is constant. Moreover, $\ap\notin\G$, because ground states and bound states require $\alpha>\ap$ by Proposition~\ref{prop:tang}(ii). Applying Lemma~\ref{lem:local}(b) at $\alpha=\ap$, we obtain $\delta_{0}$.
\end{proof}

\subsection{The pairing principle}

For $k\ge0$ we set
\[
U_{k}=\{\alpha>0:\ \N(\alpha)\ge k+1\}.
\]

\begin{proposition}\label{prop:boundary}
For every $k\ge0$, the set $U_{k}$ is open and contained in $[\ap+\delta_{0},\infty)$. Every boundary point $\hat\alpha\in(0,\infty)$ of $U_{k}$ belongs to $\G$ and satisfies $\N(\hat\alpha)\le k$.
\end{proposition}

\begin{proof}
Openness follows from Lemma~\ref{lem:local}(a), and the inclusion follows from Lemma~\ref{lem:small}. Let $\hat\alpha\in(0,\infty)$ be a boundary point of $U_{k}$. Since $U_{k}$ is open, $\hat\alpha\notin U_{k}$, that is, $\N(\hat\alpha)\le k$. If $\hat\alpha$ did not belong to $\G$, Lemma~\ref{lem:local}(b) would provide a neighbourhood of $\hat\alpha$ on which $\N\equiv\N(\hat\alpha)\le k$. This neighbourhood would not meet $U_{k}$, which contradicts the fact that $\hat\alpha$ lies in the closure of $U_{k}$.
\end{proof}

\begin{theorem}[Pairing principle]\label{thm:pairing}
Let $k\ge0$, and assume that there are $\beta_{1}<\beta_{2}$ such that $\N(\beta_{1})\ge k+1$ and $\N(\alpha)\le k$ for every $\alpha\ge\beta_{2}$. Let $(a,b)$ be the connected component of $U_{k}$ containing $\beta_{1}$. Then $\ap<a<\beta_{1}<b\le\beta_{2}$, both $a$ and $b$ belong to $\G$, $\N(a)\le k$ and $\N(b)\le k$, while $\N\ge k+1$ on $(a,b)$. In particular, \eqref{eq:P} possesses two distinct radial solutions $u(\cdot,a)$ and $u(\cdot,b)$, positive at the origin, each of which is a ground state or a bound state with at most $k$ nodes. Moreover, these two solutions are not related by a translation or by a change of sign.
\end{theorem}

\begin{proof}
By Proposition~\ref{prop:boundary}, $U_{k}$ is an open subset of $(0,\infty)$ containing $\beta_{1}$, so its connected component containing $\beta_{1}$ is an open interval $(a,b)$. The inclusion $U_{k}\subset[\ap+\delta_{0},\infty)$ gives $a\ge\ap+\delta_{0}>\ap$. Since $\N\le k$ on $[\beta_{2},\infty)$, the interval $(a,b)$ does not meet $[\beta_{2},\infty)$, and we get $b\le\beta_{2}$. Both $a$ and $b$ are finite and positive. As endpoints of a connected component of an open set, they do not belong to $U_{k}$ but lie in its closure. Hence they are boundary points of $U_{k}$, and Proposition~\ref{prop:boundary} shows that they lie in $\G$ and have at most $k$ zeros. Finally, the two solutions take the distinct positive values $a<b$ at the origin, so they are distinct radial solutions centred at the origin. A nontrivial translate of $u(\cdot,a)$ is not radially symmetric about the origin, and $-u(\cdot,a)$ is negative at the origin. We conclude that neither of them coincides with $u(\cdot,b)$.
\end{proof}

\begin{corollary}\label{cor:pairs}
Let $N_{\max}=\sup_{\alpha>0}\N(\alpha)$ and $N_{\infty}=\limsup_{\alpha\to\infty}\N(\alpha)$, and suppose that $N_{\infty}<\infty$. Then, for every integer $k$ with $N_{\infty}\le k<N_{\max}$, there are two distinct radial solutions of \eqref{eq:P}, positive at the origin, with at most $k$ nodes, each of which is a ground state or a bound state.
\end{corollary}

\begin{proof}
Since $\N$ is integer valued and $N_{\infty}\le k$, there is $\beta_{2}$ such that $\N(\alpha)\le k$ for all $\alpha\ge\beta_{2}$. Since $k<N_{\max}$, there is $\beta_{1}$ with $\N(\beta_{1})\ge k+1$, and necessarily $\beta_{1}<\beta_{2}$. Theorem~\ref{thm:pairing} applies.
\end{proof}

In order to identify exactly the number of nodes of the two solutions, we need a simple-crossing property. Lemma~\ref{lem:local}(a) already gives the corresponding lower bound unconditionally.

\begin{hypothesis}[Simple crossing]\label{hyp:S}
For every $\hat\alpha\in\G$, there is $\delta>0$ such that $\N(\beta)\le\N(\hat\alpha)+1$ for $|\beta-\hat\alpha|<\delta$.
\end{hypothesis}

\begin{corollary}\label{cor:exact}
Under the assumptions of Theorem~\ref{thm:pairing}, if Hypothesis~\ref{hyp:S} holds at $a$ and at $b$, then $\N(a)=\N(b)=k$. Consequently, if $k\ge1$, equation \eqref{eq:P} has at least two distinct $k$-node radial bound states positive at the origin, and if $k=0$, it has at least two distinct radial ground states.
\end{corollary}

\begin{proof}
Points of $(a,b)$ arbitrarily close to $a$ satisfy $\N\ge k+1$, while Hypothesis~\ref{hyp:S} gives $\N\le\N(a)+1$ near $a$. We deduce that $\N(a)\ge k$, and $\N(a)\le k$ by Theorem~\ref{thm:pairing}. The same argument applies at $b$.
\end{proof}

\begin{remark}\rm\label{rem:S}
Heuristically, Hypothesis~\ref{hyp:S} holds for the following reason. Let $\hat\alpha\in\G$ and let $\beta$ be close to $\hat\alpha$. The solution $u(\cdot,\beta)$ follows the exponentially decaying tail of $u(\cdot,\hat\alpha)$ up to a large radius $R$, at which the deviation, which grows like $e^{\sqrt\zeta r}$, becomes comparable to $e^{-\sqrt\zeta r}$. If the solution then crosses zero, its energy at the crossing is of order $e^{-2\sqrt\zeta R}$. A further zero requires an excursion to $|u|>\ap$ and back, close to the homoclinic orbit of the autonomous equation $u''+f(u)=0$. During this excursion, the energy decreases by an amount of order $R^{-1}$, which is much larger than $e^{-2\sqrt\zeta R}$. The energy therefore becomes negative, and the solution is trapped. Making this argument quantitative requires a careful perturbative analysis near the saddle point $(0,0)$, which has not been carried out.
\end{remark}

A stronger property than Hypothesis~\ref{hyp:S} gives the parity of the number of bound states. We say that $\hat\alpha\in\G$ is \emph{transversal} if there is $\delta>0$ such that $\N$ is equal to $\N(\hat\alpha)$ on one of the intervals $(\hat\alpha-\delta,\hat\alpha)$, $(\hat\alpha,\hat\alpha+\delta)$, and to $\N(\hat\alpha)+1$ on the other. A transversal point clearly satisfies Hypothesis~\ref{hyp:S}.

\begin{proposition}\label{prop:parity}
Let $k\ge0$. Assume that there is $\beta_{2}$ such that $\N(\alpha)\le k$ for all $\alpha\ge\beta_{2}$, and that every $\hat\alpha\in\G$ with $\N(\hat\alpha)\le k$ is transversal. Then every connected component of $U_{k}$ is a bounded interval $(a,b)$ with $\ap<a<b\le\beta_{2}$ and $\N(a)=\N(b)=k$. Moreover, every $\hat\alpha\in\G$ with $\N(\hat\alpha)=k$ is an endpoint of exactly one connected component of $U_{k}$. Consequently, the set of initial heights of the radial solutions of \eqref{eq:P} that are ground states (if $k=0$) or $k$-node bound states (if $k\ge1$) is either infinite or has an even number of elements.
\end{proposition}

\begin{proof}
Let $(a,b)$ be a connected component of $U_{k}$. As in the proof of Theorem~\ref{thm:pairing}, we have $\ap<a<b\le\beta_{2}$, and $a$, $b$ are points of $\G$ with $\N(a)\le k$ and $\N(b)\le k$. The point $a$ is transversal, and $\N\ge k+1>\N(a)$ on $(a,b)$. Hence $\N=\N(a)+1$ on a right neighbourhood of $a$, and we obtain $\N(a)+1\ge k+1$, that is, $\N(a)=k$. The same argument gives $\N(b)=k$.

Let now $\hat\alpha\in\G$ with $\N(\hat\alpha)=k$. By transversality, $\N=k$ on one side of $\hat\alpha$ and $\N=k+1$ on the other. Suppose, for instance, that $\N=k+1$ on $(\hat\alpha,\hat\alpha+\delta)$ and $\N=k$ on $(\hat\alpha-\delta,\hat\alpha)$. Then $(\hat\alpha,\hat\alpha+\delta)\subset U_{k}$ while $\hat\alpha\notin U_{k}$, so $\hat\alpha$ is the left endpoint of the component containing $(\hat\alpha,\hat\alpha+\delta)$. On the other hand, $\hat\alpha$ is not the right endpoint of any component, since $(\hat\alpha-\delta,\hat\alpha)$ does not meet $U_{k}$. The other case is symmetric. Therefore, the map which sends a component $(a,b)$ to the pair $\{a,b\}$ yields a partition of the set of heights in question into pairs, and the last assertion follows.
\end{proof}

\subsection{Intersection of ground states}

\begin{lemma}\label{lem:intersection}
Let $U$ and $W$ be two distinct radial ground states of \eqref{eq:P}, both centred at the origin. Then there is $r>0$ with $U(r)=W(r)$.
\end{lemma}

\begin{proof}
Two radial solutions of \eqref{eq:ivp} with the same initial height coincide, so we may assume that $U(0)<W(0)$. Suppose, by contradiction, that the graphs do not meet on $(0,\infty)$. Then $0<U<W$ on $[0,\infty)$. We consider the Wronskian $\varpi(r)=r\big(U'(r)W(r)-U(r)W'(r)\big)$. Since $(rU')'=-rf(U)$ and $(rW')'=-rf(W)$, we have
\[
\varpi'(r)=-rf(U)W+rf(W)U=rU W\left(\frac{f(W)}{W}-\frac{f(U)}{U}\right)>0
\]
for $r>0$, by Lemma~\ref{lem:elementary}(e). It follows that $\varpi$ is strictly increasing, with $\varpi(0)=0$, so that $\varpi(r)\ge\varpi(1)>0$ for all $r\ge1$. On the other hand, $U$, $U'$, $W$ and $W'$ decay exponentially by Proposition~\ref{prop:tang}(iv), hence $\varpi(r)\to0$ as $r\to\infty$. This is a contradiction.
\end{proof}

\subsection{Tang's reduction in the plane}

\begin{remark}\rm\label{rem:tangcrit}
Let $v=\partial u/\partial\alpha$. If $z(\alpha)$ is a zero of $u(\cdot,\alpha)$, the implicit function theorem shows that $z$ depends smoothly on $\alpha$, with $z'(\alpha)=-v(z,\alpha)/u'(z,\alpha)$. Tang reduces the uniqueness of bound states to two conditions on $v$ \cite[Lemma~3.1]{Tang25}. In our setting they read as follows, for a given $\alpha>\ap$.
\begin{enumerate}
\item[(a)] If $u=u(\cdot,\alpha)$ has exactly $k\ge1$ zeros $z_{1}<\dots<z_{k}$, then $v$ has exactly $k$ zeros $\tau_{1}<\dots<\tau_{k}$ in $[0,z_{k}]$, with $\tau_{1}\in(0,z_{1})$ and $\tau_{i}\in(z_{i-1},z_{i})$ for $2\le i\le k$. Moreover, $u'v'>0$ whenever $v=0$.
\item[(b)] If $u$ is a ground state, then $v$ has exactly one zero in $(0,\infty)$. If $u$ is a $k$-node bound state, then $v$ has exactly one zero in $[z_{k},\infty)$, and it lies beyond the last critical point of $u$. In either case, $v$ is strictly monotone for large $r$ and $|v(r)|\to\infty$ as $r\to\infty$.
\end{enumerate}
\end{remark}
For the power nonlinearity in dimension $n\ge3$, Tang verifies (a) and (b) in \cite[Sections~4 and~5]{Tang25}. This verification uses the structural conditions of Table~\ref{tab:conditions} and, in the transition between phases, identities that are specific to the power nonlinearity \cite[Appendix~A.4]{Tang25}. By contrast, the deduction of uniqueness from (a) and (b) is robust, as the following proposition shows. Its proof follows the argument of \cite[proof of Lemma~3.1]{Tang25}, which goes back to McLeod \cite{McLeod93}.

\begin{proposition}\label{prop:reduction}
~\begin{enumerate}
\item[(i)] Let $0<\alpha_{1}<\alpha_{2}$, and assume that condition (a) of Remark~\ref{rem:tangcrit} holds at every $\alpha\in[\alpha_{1},\alpha_{2}]$. Then $\N(\alpha_{1})\le\N(\alpha_{2})$.
\item[(ii)] Let $\hat\alpha\in\G$, and assume that conditions (a) and (b) hold at $\hat\alpha$. Then there is $\delta>0$ such that $\N(\beta)\ge\N(\hat\alpha)+1$ for every $\beta\in(\hat\alpha,\hat\alpha+\delta)$.
\item[(iii)] Assume that (a) holds at every $\alpha>\ap$ and that (b) holds at every $\alpha\in\G$. Then \eqref{eq:P} has at most one radial ground state centred at the origin and, for every $k\ge1$, at most one $k$-node radial bound state positive at the origin.
\end{enumerate}
\end{proposition}

\begin{proof}
We first show that (a) at $\alpha$ implies $u'v>0$ at every zero of $u=u(\cdot,\alpha)$. The zeros of $v$ are simple. Indeed, $v(\tau)=v'(\tau)=0$ at some $\tau>0$ would force $v\equiv0$ by uniqueness for the linear equation of $v$, which contradicts $v(0)=1$. Hence $v$ changes sign at each $\tau_{i}$, and $v$ has the sign $(-1)^{i}$ on $(\tau_{i},\tau_{i+1})\cap[0,z_{k}]$, where $\tau_{0}=0$ and $\tau_{k+1}=\infty$. Since $z_{i}\in(\tau_{i},\tau_{i+1})$, we get that $v(z_{i})$ has the sign $(-1)^{i}$. On the other hand, $u>0$ on $[0,z_{1})$, the zeros of $u$ are simple, and therefore $u'(z_{i})$ also has the sign $(-1)^{i}$. We conclude that $u'v>0$ at each $z_{i}$.

(i) We may assume that $m=\N(\alpha_{1})\ge1$. Let $z_{1}<\dots<z_{m}$ be the zeros of $u(\cdot,\alpha_{1})$. By continuity of $u$ on $[0,1]\times[\alpha_{1},\alpha_{2}]$ and $u(0,\alpha)=\alpha\ge\alpha_{1}$, there is $r_{0}>0$ such that $u(r,\alpha)>0$ for $r\in[0,r_{0}]$ and $\alpha\in[\alpha_{1},\alpha_{2}]$. By the implicit function theorem, each $z_{j}$ is the starting point of a $C^{1}$ curve $\alpha\mapsto\zeta_{j}(\alpha)$ of zeros of $u(\cdot,\alpha)$, with $\zeta_{j}'=-v/u'<0$ by the first step. Let $\alpha^{\sharp}$ be the supremum of those $\beta\in[\alpha_{1},\alpha_{2}]$ such that the $m$ curves are defined on $[\alpha_{1},\beta]$. Each $\zeta_{j}$ is decreasing and takes values in $[r_{0},z_{j}]$, so it has a limit as $\alpha\uparrow\alpha^{\sharp}$, which is a zero of $u(\cdot,\alpha^{\sharp})$ by continuity. These limits are distinct. Indeed, near a simple zero of $u(\cdot,\alpha^{\sharp})$ the zero set of $u$ is a single $C^{1}$ curve, so two of the curves cannot converge to the same point. If $\alpha^{\sharp}<\alpha_{2}$, the implicit function theorem extends the curves beyond $\alpha^{\sharp}$, which contradicts the definition of $\alpha^{\sharp}$. Hence the curves are defined on $[\alpha_{1},\alpha_{2}]$, and $u(\cdot,\alpha_{2})$ has at least $m$ zeros.

(ii) Let $j=\N(\hat\alpha)$, $\hat u=u(\cdot,\hat\alpha)$ and $\hat v=v(\cdot,\hat\alpha)$. We assume that $j$ is odd; the case of even $j$, including $j=0$, follows by reflection. Then $\hat u<0$ eventually, and $\hat u$ increases to $0$ beyond its last critical point. By (a) and (b), $\hat v$ has exactly $j+1$ zeros, all simple, so that $\hat v>0$ eventually. Since $\hat v$ is eventually monotone and $|\hat v|\to\infty$, we deduce that $\hat v,\hat v'>0$ for large $r$. Let $s_{0}>0$ be such that $f'<0$ on $(-s_{0},s_{0})$, which is possible since $f'(0)<0$. We fix $\bar r$ beyond the last zero and the last critical point of $\hat u$, such that $\hat v(\bar r)>0$, $\hat v'(\bar r)>0$ and $|\hat u|<s_{0}$ on $[\bar r,\infty)$. Let $\beta>\hat\alpha$, let $u_{\beta}=u(\cdot,\beta)$ and $w=u_{\beta}-\hat u$. Since $w(\bar r)=\int_{\hat\alpha}^{\beta}v(\bar r,s)\dd s$ and $w'(\bar r)=\int_{\hat\alpha}^{\beta}v'(\bar r,s)\dd s$, the continuity of $v$ and $v'$ gives $w(\bar r)>0$ and $w'(\bar r)>0$ for $\beta$ close to $\hat\alpha$. Arguing as in the proof of Lemma~\ref{lem:local}(b), we also see that $u_{\beta}$ has exactly $j$ zeros in $(0,\bar r)$ and $u_{\beta}(\bar r)<0$ for such $\beta$. Suppose that $u_{\beta}$ has no zero in $(\bar r,\infty)$. If $\beta\in\G$, then $w(r)\to0$ as $r\to\infty$. Otherwise, $u_{\beta}$ is oscillatory and negative on $(\bar r,\infty)$, so it oscillates about $-u_{1}$ by Proposition~\ref{prop:tang}(iv), and it takes values below $-u_{1}$ at arbitrarily large radii, while $\hat u(r)\to0$. In this case, $w$ takes negative values at arbitrarily large radii. In both cases, since $w(\bar r)>0$ and $w'(\bar r)>0$, the function $w$ attains a positive local maximum at some $\hat r>\bar r$, where $w'(\hat r)=0$ and $w''(\hat r)\le0$. Now $w$ solves
\[
w''+\frac{1}{r}w'+c(r)w=0,\qquad c(r)=\int_{0}^{1}f'\big(\hat u(r)+\theta w(r)\big)\dd\theta .
\]
At $\hat r$ we have $\hat u(\hat r)<\hat u(\hat r)+\theta w(\hat r)\le u_{\beta}(\hat r)<0$ and $|\hat u(\hat r)|<s_{0}$, so that $c(\hat r)<0$. We infer that $w''(\hat r)=-c(\hat r)w(\hat r)>0$, a contradiction. Hence $u_{\beta}$ has a zero in $(\bar r,\infty)$, and $\N(\beta)\ge j+1$.

(iii) Suppose that $\hat\alpha<\hat\alpha'$ are two points of $\G$ with $\N(\hat\alpha)=\N(\hat\alpha')=k$. By (ii), there is $\beta\in(\hat\alpha,\hat\alpha')$ with $\N(\beta)\ge k+1$. By (i), applied on $[\beta,\hat\alpha']\subset(\ap,\infty)$, we get $\N(\hat\alpha')\ge\N(\beta)\ge k+1$, a contradiction.
\end{proof}

\begin{corollary}\label{cor:failure}
Under the assumptions of Theorem~\ref{thm:pairing}, let $(a,b)$ be the interval given there. Then, for every $\delta\in(0,b-a)$, condition (a) of Remark~\ref{rem:tangcrit} fails at some $\alpha\in[b-\delta,b]$.
\end{corollary}

\begin{proof}
We have $\N(b-\delta)\ge k+1>k\ge\N(b)$, and Proposition~\ref{prop:reduction}(i) applies on $[b-\delta,b]$.
\end{proof}

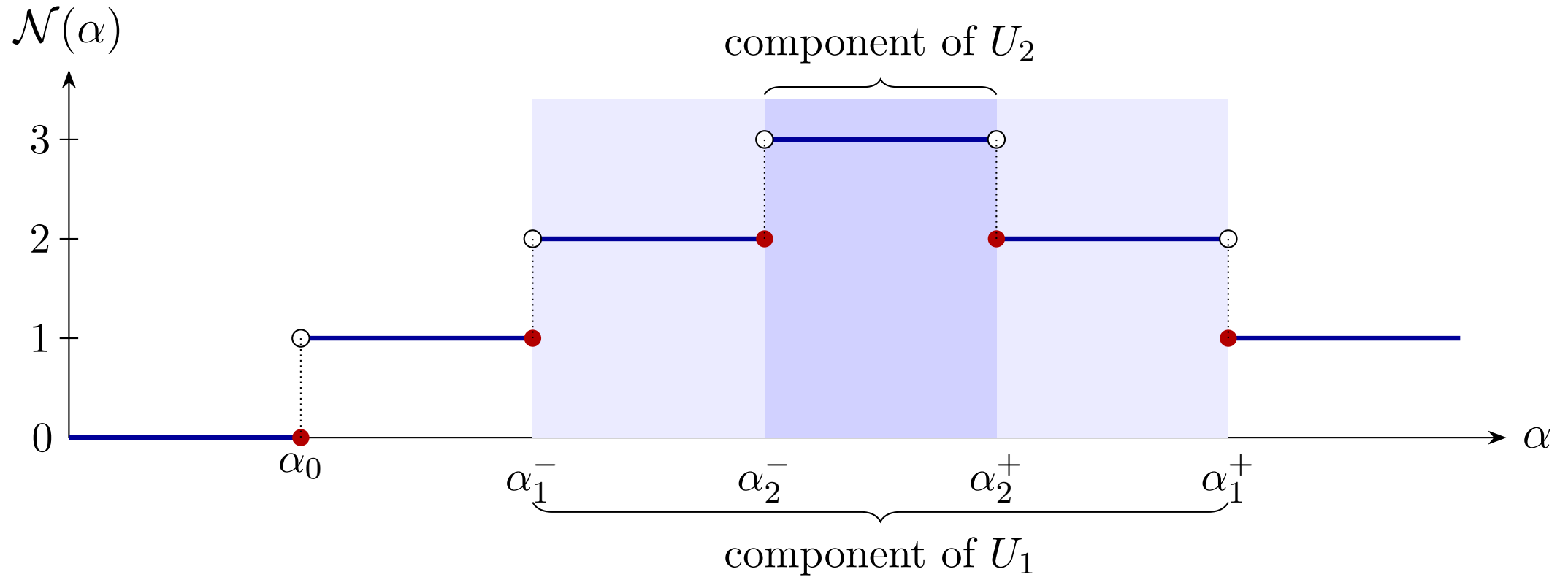
\begin{figure}[!htbp]
\centering
\begin{tikzpicture}[x=1.05cm,y=0.9cm,>=Stealth]
\draw[->] (0,0)--(12.4,0) node[right] {$\alpha$};
\draw[->] (0,0)--(0,3.7) node[above] {$\N(\alpha)$};
\foreach \y in {1,2,3} \draw (-0.08,\y)--(0.08,\y) node[left=3pt] {\small $\y$};
\node[left] at (0,0) {\small $0$};
\fill[blue!8] (4,0) rectangle (10,3.4);
\fill[blue!18] (6,0) rectangle (8,3.4);
\draw[very thick,blue!60!black] (0,0)--(2,0);
\draw[very thick,blue!60!black] (2,1)--(4,1);
\draw[very thick,blue!60!black] (4,2)--(6,2);
\draw[very thick,blue!60!black] (6,3)--(8,3);
\draw[very thick,blue!60!black] (8,2)--(10,2);
\draw[very thick,blue!60!black] (10,1)--(12,1);
\foreach \p in {(2,0),(4,1),(6,2),(8,2),(10,1)} \fill[red!70!black] \p circle (2.2pt);
\foreach \p in {(2,1),(4,2),(6,3),(8,3),(10,2)} \draw[fill=white] \p circle (2.2pt);
\draw[densely dotted] (2,0)--(2,1) (4,1)--(4,2) (6,2)--(6,3) (8,2)--(8,3) (10,1)--(10,2);
\node[below] at (2,0) {\small $\alpha_{0}$};
\node[below] at (4,0) {\small $\alpha_{1}^{-}$};
\node[below] at (6,0) {\small $\alpha_{2}^{-}$};
\node[below] at (8,0) {\small $\alpha_{2}^{+}$};
\node[below] at (10,0) {\small $\alpha_{1}^{+}$};
\draw[decorate,decoration={brace,amplitude=5pt,mirror}] (4,-0.65)--(10,-0.65) node[midway,below=5pt] {\small component of $U_{1}$};
\draw[decorate,decoration={brace,amplitude=4pt}] (6,3.45)--(8,3.45) node[midway,above=4pt] {\small component of $U_{2}$};
\end{tikzpicture}
\caption{Schematic graph of $\N$ for $\lambda$ close to $1$, as predicted by Conjectures~\ref{conj:large} and~\ref{conj:BL}; it is not computed. Filled dots are the values of $\N$ at points of $\G$. The endpoints of each component of $U_{k}$ are radial solutions with at most $k$ nodes (Theorem~\ref{thm:pairing}), here $k$-node bound states, so that bound states come in nested pairs $\alpha_{k}^{-}<\alpha_{k}^{+}$. Tang's conditions force $\N$ to be non-decreasing, and therefore they must fail near $\alpha_{1}^{+}$ and $\alpha_{2}^{+}$ (Corollary~\ref{cor:failure}).}\label{fig:pairing}
\end{figure}

In particular, if \eqref{eq:P} has two distinct $k$-node bound states for some $k\ge1$, then Tang's conditions cannot hold for all initial heights. The mechanism is illustrated in Figure~\ref{fig:pairing}.

\section{Bubbling as the initial height tends to infinity}\label{sec:bubble}

In this section, we describe the behaviour of $u(\cdot,\alpha)$ as $\alpha\to\infty$. Lemmas~\ref{lem:liouville} and~\ref{lem:core} are rigorous, while the matched asymptotics built on them in Subsections~\ref{sub:z1} and~\ref{sub:mechanism} are formal and are labelled as such.

\subsection{The rescaled equation}
Let $\alpha>0$, and define the scale $\varepsilon_{\alpha}$ and the rescaled function $w_{\alpha}$ by
\[
\varepsilon_{\alpha}=\big(2\lambda\alpha^{2}e^{\alpha^{2}}\big)^{-1/2},\qquad u(\varepsilon_{\alpha}s,\alpha)=\alpha+\frac{w_{\alpha}(s)}{2\alpha}.
\]
Then $u^{2}=\alpha^{2}+w_{\alpha}+w_{\alpha}^{2}/(4\alpha^{2})$. Since $2\alpha^{2}\varepsilon_{\alpha}^{2}=\lambda^{-1}e^{-\alpha^{2}}$, a direct computation shows that \eqref{eq:ivp} becomes
\begin{equation}\label{eq:rescaled}
-\Delta_{s}w_{\alpha}=\Big(1+\frac{w_{\alpha}}{2\alpha^{2}}\Big)\Big(e^{\,w_{\alpha}+w_{\alpha}^{2}/(4\alpha^{2})}-\frac{e^{-\alpha^{2}}}{\lambda}\Big),\qquad w_{\alpha}(0)=0,\ w_{\alpha}'(0)=0.
\end{equation}
The linear term $-u$ of the original equation contributes only the term of order $e^{-\alpha^{2}}$. Expanding the right-hand side of \eqref{eq:rescaled} in powers of $\alpha^{-2}$, we are led to the expansion $w_{\alpha}=w_{0}+\alpha^{-2}w_{1}+o(\alpha^{-2})$, which is justified on compact sets in Lemma~\ref{lem:core} below. Here
\[
-\Delta w_{0}=e^{w_{0}},\qquad -\Delta w_{1}-e^{w_{0}}w_{1}=e^{w_{0}}\Big(\frac{w_{0}^{2}}{4}+\frac{w_{0}}{2}\Big),
\]
with zero initial data at $s=0$. The first equation gives $w_{0}(s)=-2\log(1+s^{2}/8)$. This is the radial solution of the Liouville equation $-\Delta w=e^{w}$ with $w(0)=0$, and it satisfies $\int_{\R^{2}}e^{w_{0}}=8\pi$. All solutions of the Liouville equation with $e^{w}\in L^{1}(\R^{2})$ were classified by Chen and Li \cite{ChenLi91}. We refer to Brezis and Merle \cite{BrezisMerle91} for the general blow-up theory of Liouville-type equations, and to \cite{AdimurthiStruwe00,Druet06,MalchiodiMartinazzi14,MancMart17,GrossiNaimen20} for bubbling at Trudinger--Moser criticality.

\subsection{The second-order correction}
The asymptotic flux $\lim_{s\to\infty}sw'(s)$ of the profile governs the behaviour of $u$ in the neck region. For $w_{0}$ it equals $-4$, by the explicit formula. The next lemma computes exactly the correction of order $\alpha^{-2}$.

\begin{lemma}\label{lem:liouville}
Let $w_{0}(s)=-2\log(1+s^{2}/8)$ and let $w_{1}$ be the solution of
\begin{equation}\label{eq:w1}
w_{1}''+\frac{1}{s}w_{1}'+e^{w_{0}}w_{1}=-e^{w_{0}}\Big(\frac{w_{0}^{2}}{4}+\frac{w_{0}}{2}\Big),\quad s>0,\qquad w_{1}(0)=w_{1}'(0)=0.
\end{equation}
Then $w_{1}(s)=O(\log s)$ as $s\to\infty$, and $\lim_{s\to\infty}s\,w_{1}'(s)=-4$.
\end{lemma}

\begin{proof}
We set $g=e^{w_{0}}(w_{0}^{2}/4+w_{0}/2)$ and $t=s^{2}/8$, so that $e^{w_{0}}=(1+t)^{-2}$ and $g=(1+t)^{-2}\big(\ell^{2}-\ell\big)$ with $\ell=\log(1+t)$. Equation \eqref{eq:w1} is linear with smooth coefficients, and its regular singular point at $s=0$ is handled exactly as in \eqref{eq:ivp}. Hence it has a unique global solution. In the variable $\tau=\log s$, the function $y(\tau)=w_{1}(e^{\tau})$ satisfies $y''=-e^{2\tau}\big(e^{w_{0}}y+g\big)$. For $\tau\ge0$ we have $e^{2\tau}e^{w_{0}}\le64e^{-2\tau}$ and $e^{2\tau}|g|\le C(1+\tau^{2})e^{-2\tau}$. Integrating twice, we obtain for $\tau\ge0$
\[
|y(\tau)|\le|y(0)|+\tau|y'(0)|+\int_{0}^{\tau}(\tau-\sigma)\Big(64e^{-2\sigma}|y(\sigma)|+C(1+\sigma^{2})e^{-2\sigma}\Big)\dd\sigma .
\]
Since $(\tau-\sigma)/(1+\tau)\le1$ for $0\le\sigma\le\tau$, the function $q(\tau)=|y(\tau)|/(1+\tau)$ satisfies $q(\tau)\le C'+\int_{0}^{\tau}64(1+\sigma)e^{-2\sigma}q(\sigma)\dd\sigma$. By Gronwall's lemma, $q$ is bounded, that is, $w_{1}(s)=O(\log s)$. We infer that $y''$ is integrable on $(0,\infty)$, and therefore $m_{1}=\lim_{s\to\infty}sw_{1}'(s)=\lim_{\tau\to\infty}y'(\tau)$ exists.

To compute $m_{1}$, we use the function $\varphi_{0}=2+sw_{0}'(s)=2(8-s^{2})/(8+s^{2})$. It is the derivative at $\mu=1$ of the family $w_{0}(\mu s)+2\log\mu$ of solutions of the Liouville equation, hence $\varphi_{0}''+\varphi_{0}'/s+e^{w_{0}}\varphi_{0}=0$. We multiply \eqref{eq:w1} by $s\varphi_{0}$, subtract the equation of $\varphi_{0}$ multiplied by $sw_{1}$, and integrate over $(0,S)$. We get
\[
S\big(\varphi_{0}(S)w_{1}'(S)-w_{1}(S)\varphi_{0}'(S)\big)=-\int_{0}^{S}s\varphi_{0}(s)g(s)\dd s .
\]
As $S\to\infty$, we have $\varphi_{0}(S)\to-2$, $S\varphi_{0}'(S)=O(S^{-2})$ and $w_{1}(S)=O(\log S)$. Passing to the limit, we deduce that $-2m_{1}=-\int_{0}^{\infty}s\varphi_{0}g\dd s$. In the variable $t$, we have $s\dd s=4\dd t$ and $\varphi_{0}=2(1-t)/(1+t)$, so that
\[
m_{1}=\frac12\int_{0}^{\infty}s\varphi_{0}g\dd s=4\int_{0}^{\infty}\frac{1-t}{(1+t)^{3}}\big(\ell^{2}-\ell\big)\dd t=4\int_{1}^{\infty}\frac{2-x}{x^{3}}\big(\log^{2}x-\log x\big)\dd x .
\]
Using $\int_{1}^{\infty}x^{-k}\log^{j}x\dd x=j!/(k-1)^{j+1}$ for $k>1$, we find
\[
\int_{1}^{\infty}\frac{2-x}{x^{3}}\log^{2}x\dd x=\frac12-2=-\frac32,\qquad \int_{1}^{\infty}\frac{2-x}{x^{3}}\log x\dd x=\frac12-1=-\frac12 .
\]
We conclude that $m_{1}=4\big(-\frac32+\frac12\big)=-4$.
\end{proof}

The expansion of $w_{\alpha}$ is rigorous on compact sets of the core variable.

\begin{lemma}\label{lem:core}
Let $S>0$. There is $\alpha_{S}>0$ such that, for every $\alpha\ge\alpha_{S}$, the solution $u(\cdot,\alpha)$ is positive on $[0,S\varepsilon_{\alpha}]$. Moreover, as $\alpha\to\infty$,
\[
w_{\alpha}=w_{0}+\alpha^{-2}w_{1}+o(\alpha^{-2})\qquad\text{in }C^{1}([0,S]).
\]
In particular, $z_{1}(\alpha)/\varepsilon_{\alpha}\to\infty$ as $\alpha\to\infty$ for the first zero $z_{1}(\alpha)$ of $u(\cdot,\alpha)$, whenever it exists.
\end{lemma}

\begin{proof}
We set $\eta=\alpha^{-2}$ and, for $w\in\R$ and $\eta\ge0$,
\[
G(w,\eta)=\Big(1+\frac{\eta w}{2}\Big)\Big(e^{\,w+\eta w^{2}/4}-\frac{e^{-1/\eta}}{\lambda}\Big)\quad(\eta>0),\qquad G(w,0)=e^{w}.
\]
For $\eta<0$ we set $G(w,\eta)=(1+\eta w/2)e^{\,w+\eta w^{2}/4}$. Since $\eta\mapsto e^{-1/\eta}$, extended by zero for $\eta\le0$, is smooth on $\R$ with all derivatives vanishing at $0$, the function $G$ is of class $C^{1}$ on $\R^{2}$. Moreover, $G_{w}(w,0)=e^{w}$ and $G_{\eta}(w,0)=e^{w}(w/2+w^{2}/4)$. Let $X=C([0,S])$ and let $\mathcal{K}$ be the Volterra operator $(\mathcal{K}h)(s)=\int_{0}^{s}\tau\log(s/\tau)h(\tau)\dd\tau$, whose kernel is bounded on $\{0\le\tau\le s\le S\}$. By \eqref{eq:rescaled}, $w_{\alpha}$ is the unique solution in $X$ of $\Phi(W,\alpha^{-2})=0$, where
\[
\Phi(W,\eta)=W+\mathcal{K}\big(G(W,\eta)\big),\qquad (W,\eta)\in X\times\R .
\]
Uniqueness holds because $G$ is locally Lipschitz and $\mathcal{K}$ has a bounded kernel, so that Gronwall's lemma applies. The map $\Phi$ is of class $C^{1}$, since $G$ is, and $\Phi(w_{0},0)=0$. Its partial differential $D_{W}\Phi(w_{0},0)=I+\mathcal{K}(e^{w_{0}}\,\cdot\,)$ is the identity plus a Volterra operator with bounded kernel, hence it is invertible on $X$. By the implicit function theorem, there are $\eta_{0}>0$ and a $C^{1}$ map $\eta\mapsto W(\cdot,\eta)\in X$ on $(-\eta_{0},\eta_{0})$ with $W(\cdot,0)=w_{0}$ and $\Phi(W(\cdot,\eta),\eta)=0$. By uniqueness, $w_{\alpha}=W(\cdot,\alpha^{-2})$ on $[0,S]$ for $\alpha^{-2}<\eta_{0}$. Differentiating $\Phi(W(\cdot,\eta),\eta)=0$ at $\eta=0$, we find that $Y=\partial_{\eta}W(\cdot,0)$ solves $Y+\mathcal{K}\big(e^{w_{0}}Y+e^{w_{0}}(w_{0}/2+w_{0}^{2}/4)\big)=0$, which is the integral form of \eqref{eq:w1}. Hence $Y=w_{1}$, and $w_{\alpha}=w_{0}+\alpha^{-2}w_{1}+o(\alpha^{-2})$ in $X$. Finally, $W'(s)=-s^{-1}\int_{0}^{s}\tau G(W(\tau),\eta)\dd\tau$, and the map $h\mapsto s^{-1}\int_{0}^{s}\tau h\dd\tau$ is bounded on $X$. We conclude that the expansion also holds in $C^{1}([0,S])$. In particular, $|w_{\alpha}|\le C_{S}$ on $[0,S]$ for $\alpha$ large, so that $u(\varepsilon_{\alpha}s,\alpha)\ge\alpha-C_{S}/(2\alpha)>0$ for $s\in[0,S]$ and $\alpha$ large. Since $S$ is arbitrary, it follows that $z_{1}(\alpha)/\varepsilon_{\alpha}\to\infty$.
\end{proof}

\subsection{Formal asymptotics of the first zero}\label{sub:z1}
The arguments of this subsection and of the next one are formal. By Lemmas~\ref{lem:liouville} and~\ref{lem:core}, the expected asymptotic flux of $w_{\alpha}$ is $-4-4\alpha^{-2}+O(\alpha^{-4})$. In terms of $u$, this means that $ru'(r)=-\frac{2}{\alpha}\big(1+\alpha^{-2}+O(\alpha^{-4})\big)$ in the neck region $\varepsilon_{\alpha}\ll r\ll\alpha^{-1}$. In this region, the right-hand side of \eqref{eq:rescaled}, multiplied by $s^{2}$, is of order $\exp\{-2\log s\,(1-2\log s/\alpha^{2})\}$. It is therefore exponentially small in the interior of the neck, whose length in the variable $\log r$ is close to $\alpha^{2}/2$, so that the flux is essentially conserved there. The function $u$ is then almost harmonic in $\log r$. We use $\log(1/\varepsilon_{\alpha})=\frac{\alpha^{2}}{2}+\frac12\log(2\lambda\alpha^{2})$ together with $w_{0}(s)=-4\log s+2\log8+o(1)$, and we obtain
\[
u(r,\alpha)=\frac{1}{\alpha}\Big(-2\log r-\log(2\lambda\alpha^{2})+\log8-1\Big)+o\big(\alpha^{-1}\big)
\]
up to $r$ of order $\alpha^{-1}$. The constant $-1$ comes from the product of the correction $-4\alpha^{-2}$ of the flux with the length $\alpha^{2}/2$ of the neck. Without this correction, one would obtain the leading-order prediction $z_{1}\approx\sqrt8/(\alpha\sqrt{2\lambda})$. The linear term of the equation modifies the flux by a relative amount $O(\alpha^{-2}\log\alpha)$ only, over a range of $\log r$ of length $O(1)$. Taking the correction into account, we are led to
\begin{equation}\label{eq:z1}
\lim_{\alpha\to\infty}\alpha\sqrt{2\lambda}\,z_{1}(\alpha)=\sqrt8\,e^{-1/2},
\end{equation}
independently of $\lambda\in(0,1)$.

\subsection{A formal mechanism for a bounded zero count}\label{sub:mechanism}
Beyond the first zero, in the range $z_{1}\ll r\ll1$, the formal expansion above gives $u(r)\approx-\frac{2}{\alpha}\log(r/z_{1})$ and $u'(r)\approx-\frac{2}{\alpha r}$, while $F(u)=-\frac{\zeta}{2}u^{2}+O(u^{4})$. Hence
\[
E(r)\approx\frac{2}{\alpha^{2}}\Big(\frac{1}{r^{2}}-\zeta\log^{2}\frac{r}{z_{1}}\Big).
\]
The energy therefore vanishes at a radius $r_{\alpha}$ with $r_{\alpha}\sqrt\zeta\log(r_{\alpha}/z_{1})=1$, that is, $r_{\alpha}\sim(\sqrt\zeta\log\alpha)^{-1}$. This is consistent with the approximations used, provided that $\sqrt\zeta\log\alpha\gg1$. Indeed, $|u(r_{\alpha})|$ is then of order $\alpha^{-1}\log\alpha$, and the contribution of the linear term to the flux on $(z_{1},r_{\alpha})$ is of relative order $(\log\alpha)^{-1}$. Since $E$ is decreasing and positive at every zero, the solution has no zero beyond $r_{\alpha}$. This suggests that $\N(\alpha)=1$ for all sufficiently large $\alpha$, for every $\lambda\in(0,1)$. It also suggests that the large-height regime sets in only when $\log\alpha\gg(1-\lambda)^{-1/2}$, which leaves room for many oscillations at moderate heights when $\lambda$ is close to $1$.

\begin{figure}[!htbp]
\centering
\begin{tikzpicture}[x=1.15cm,y=0.85cm,>=Stealth]
\draw[->] (-0.2,0)--(11.6,0) node[right] {$\log r$};
\draw[->] (0,-1.3)--(0,4.2) node[above] {$u(r,\alpha)$};
\draw[dashed,gray] (0,3.4)--(1.8,3.4);
\node[left] at (0,3.4) {\small $\alpha$};
\draw[very thick,blue!60!black,domain=0:11.3,samples=150] plot (\x,{3.4-0.78*(max(2*(\x-2.4),0)+ln(1+exp(-abs(2*(\x-2.4)))))/2+0.6*(max(1.5*(\x-7.6),0)+ln(1+exp(-abs(1.5*(\x-7.6)))))/1.5});
\foreach \x/\t in {2.2/{\varepsilon_{\alpha}},6.4/{\alpha^{-1}},9.3/{(\log\alpha)^{-1}}} {\draw[densely dotted] (\x,-1.3)--(\x,4.0); \node[below right] at (\x,-1.3) {\small $\t$};}
\fill[red!70!black] (6.87,0) circle (2.2pt) node[above right] {\small $z_{1}$};
\node at (1.1,4.0) {\small core};
\node at (4.3,4.0) {\small neck};
\node at (7.85,4.0) {\small transition};
\node at (10.45,4.0) {\small $E<0$};
\node[align=left,font=\footnotesize] at (3.7,0.9) {$ru'\approx-\tfrac{2}{\alpha}(1+\alpha^{-2})$};
\end{tikzpicture}
\caption{Schematic structure of $u(\cdot,\alpha)$ for large $\alpha$, plotted against $\log r$ (not to scale). In the core $r\lesssim\varepsilon_{\alpha}$, the solution is a Liouville bubble (Lemmas~\ref{lem:liouville} and~\ref{lem:core}, rigorous). In the neck, $u$ is almost linear in $\log r$ with the slope indicated, and it vanishes for the first time near $z_{1}\approx\sqrt8\,e^{-1/2}/(\alpha\sqrt{2\lambda})$. Beyond $r_{\alpha}\sim(\sqrt\zeta\log\alpha)^{-1}$, the energy is negative, no further zero can occur, and the solution eventually oscillates about $-u_{1}$ (Proposition~\ref{prop:tang}(iv)). The neck and transition regions are described only formally.}\label{fig:regions}
\end{figure}
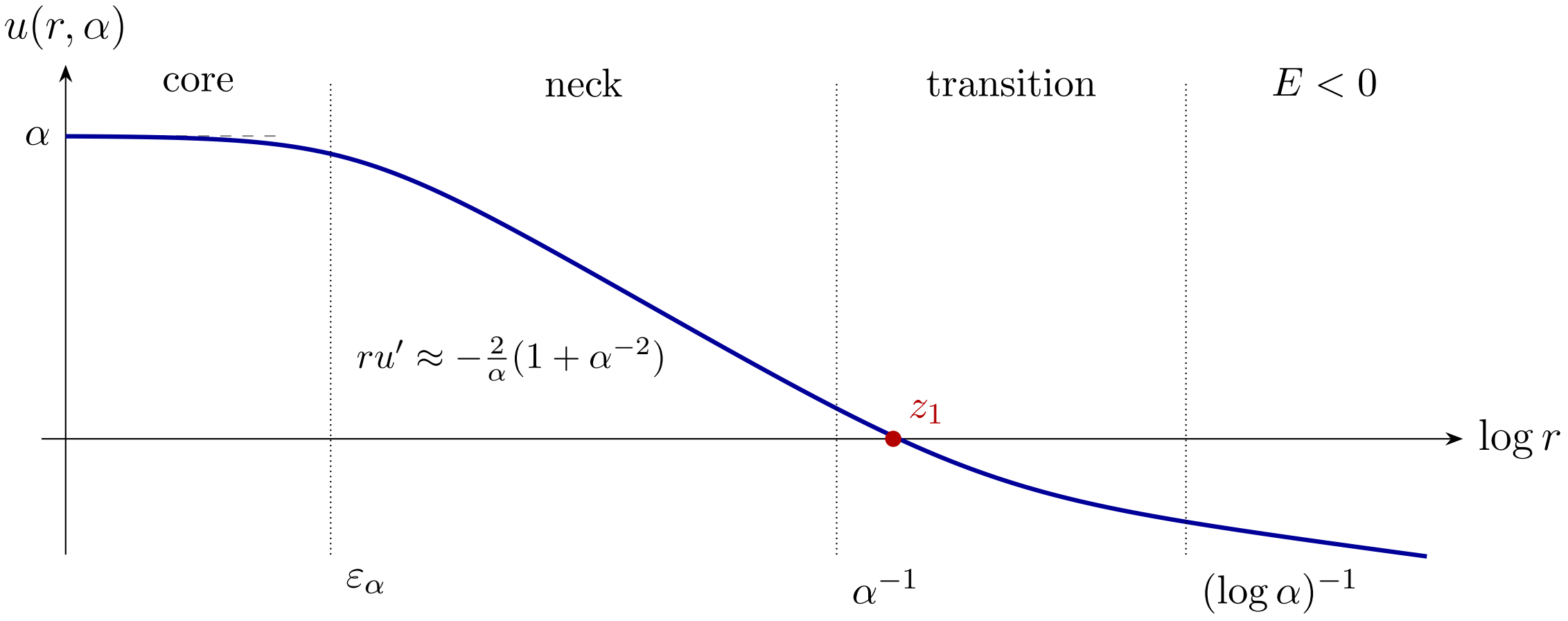

The three regions involved in this description are sketched in Figure~\ref{fig:regions}.

\section{Conjectures}\label{sec:conjectures}

The following statements are not proved. They are motivated by the results of Sections~\ref{sec:shooting} and~\ref{sec:bubble}, by the analogy with the Dirichlet problem in a disc discussed in the introduction and, for Conjectures~\ref{conj:threshold} and~\ref{conj:BL}, by numerical experiments that are not reported here.

\begin{conjecture}[Ground state]\label{conj:ground}
For every $\lambda\in(0,1)$, equation \eqref{eq:P} has exactly one positive solution up to translation, and this solution is nondegenerate: the kernel in $L^{2}(\R^{2})$ of the linearised operator $\Delta+f'(u)=\Delta-1+\lambda e^{u^{2}}(1+2u^{2})$ is spanned by $\partial_{1}u$ and $\partial_{2}u$.
\end{conjecture}

\begin{conjecture}[Large heights]\label{conj:large}
For every $\lambda\in(0,1)$, one has $\N(\alpha)=1$ for all sufficiently large $\alpha$. Moreover, \eqref{eq:z1} holds.
\end{conjecture}

\begin{conjecture}[Threshold]\label{conj:threshold}
There is $\lambda_{c}\in(0,1)$ such that \eqref{eq:P} has no radial sign-changing solution in $H^{1}(\R^{2})$ when $0<\lambda\le\lambda_{c}$, and has radial bound states when $\lambda_{c}<\lambda<1$. Moreover, the maximal number of nodes of a radial bound state tends to infinity as $\lambda\to1$.
\end{conjecture}

\begin{conjecture}[Failure of the Berestycki--Lions uniqueness]\label{conj:BL}
Let $\lambda_{c}<\lambda<1$ and $k\ge1$. The number of $k$-node radial bound states of \eqref{eq:P} that are positive at the origin is even (possibly zero), and all of them are nondegenerate. More precisely, their initial heights form nested pairs $\alpha_{k}^{-}<\alpha_{k}^{+}$: the zero count $\N$ jumps from $k$ to $k+1$ at $\alpha_{k}^{-}$ and from $k+1$ to $k$ at $\alpha_{k}^{+}$. Tang's conditions (a) and (b) of Remark~\ref{rem:tangcrit} hold at $\alpha_{k}^{-}$. In particular, the uniqueness of $k$-node bound states fails whenever they exist.
\end{conjecture}

Several parts of Conjecture~\ref{conj:BL} reduce to the other conjectures and to the rigorous results of Section~\ref{sec:shooting}. First, Conjecture~\ref{conj:large} gives $N_{\infty}=1\le k$. Hence, by Corollaries~\ref{cor:pairs} and~\ref{cor:exact}, the existence of two $k$-node bound states would follow from Conjecture~\ref{conj:large}, Hypothesis~\ref{hyp:S} and the existence of a height $\beta_{1}$ with $\N(\beta_{1})\ge k+1$. Second, by Proposition~\ref{prop:parity}, the parity statement would follow from Conjecture~\ref{conj:large} and the transversality of every point $\hat\alpha\in\G$ with $\N(\hat\alpha)\le k$. We expect transversality to follow from the nondegeneracy of all bound states, but this implication has not been proved. Third, if the structure described in Conjecture~\ref{conj:BL} holds, then $\N=k$ on a right neighbourhood of $\alpha_{k}^{+}$. By Proposition~\ref{prop:reduction}(ii), conditions (a) and (b) cannot both hold at $\alpha_{k}^{+}$.

\section{Towards proofs}\label{sec:program}

\subsection{Large heights}
The core of the programme is Conjecture~\ref{conj:large}. The tools are the blow-up analysis for Trudinger--Moser critical equations \cite{AdimurthiStruwe00,Druet06,MalchiodiMartinazzi14,MancMart17,GrossiNaimen20}, the classification of the limiting Liouville profile \cite{ChenLi91} and the quantisation of the blow-up mass \cite{BrezisMerle91,Druet06}. One has to establish a precise expansion of $u(\cdot,\alpha)$ in three regions. The first is the core $r\lesssim\varepsilon_{\alpha}$, where Lemma~\ref{lem:core} already gives the expansion on compact sets of the rescaled variable. The second is the neck $\varepsilon_{\alpha}\ll r\lesssim\alpha^{-1}$, where the solution is almost harmonic in $\log r$. The third is the region $\alpha^{-1}\lesssim r\lesssim(\log\alpha)^{-1}$, where the solution is small and the energy changes sign. The expansion must be accurate enough to show that $E(r_{\alpha})<0$ at a radius $r_{\alpha}$ beyond the first zero, as in the formal argument of Section~\ref{sec:bubble}. Lemma~\ref{lem:liouville} identifies the second-order term of the core expansion, and the factor $e^{-1/2}$ in \eqref{eq:z1} is a useful test of any such expansion. Once Conjecture~\ref{conj:large} is known, Corollary~\ref{cor:pairs} produces, for every $\lambda$ with $N_{\max}\ge2$, two distinct radial solutions with at most one node.

\subsection{Finite heights and validated numerics}
The remaining hypothesis of Theorem~\ref{thm:pairing}, namely the existence of $\beta_{1}$ with $\N(\beta_{1})\ge k+1$, concerns a single solution on a compact interval. It is therefore accessible to rigorous computer-assisted proofs in the spirit of Cohen, Li and Schlag \cite{CLS24}. Indeed, it suffices to enclose $u(\cdot,\beta_{1})$ on an interval $[0,R]$ containing $k+1$ zeros. Combined with a proof of Conjecture~\ref{conj:large}, this would give, for the corresponding value of $\lambda$, a rigorous proof that \eqref{eq:P} has two distinct radial solutions with at most $k$ nodes. Under Hypothesis~\ref{hyp:S}, these would be two distinct $k$-node bound states.

\subsection{Uniqueness of the ground state}
Lemma~\ref{lem:obstructions} rules out a direct use of the Serrin--Tang criterion and of Tang's first-phase lemma, and two features of the problem suggest a hybrid strategy. In the core, where the exponential growth dominates, the change of variables $r=2e^{-t/2}$ transforms the radial equation $u''+u'/r+f(u)=0$ into $u_{tt}+e^{-t}f(u)=0$. This is the Emden--Fowler form used by Atkinson and Peletier \cite{AtkinsonPeletier86,AtkinsonPeletier86b}, on which the uniqueness proofs of Adimurthi \cite{Adimurthi98} and Tarsi \cite{Tarsi03} for the Dirichlet problem in a disc are based. Lemma~\ref{lem:intersection} is the whole-plane version of the first step of an intersection argument. A natural next step would be to show that two ground states intersect exactly once, and that a suitable comparison quantity is monotone beyond the intersection point, in the spirit of the separation method of Peletier and Serrin \cite{PeletierSerrin83}. In the tail, where $f(s)=-(1-\lambda)s+O(s^{3})$, the variation analysis of \cite{Tang25} is a natural candidate. Matching the two regimes is the main difficulty, and we have not attempted it.

\subsection{Further questions}
The same questions arise for general nonlinearities with critical exponential growth satisfying the hypotheses of \cite{dFMR95,AMS12,RufSani13}, for the nonlinearity $(e^{4\pi u^{2}}-1-4\pi u^{2})u$ of \cite{DKM20} (see Remark~\ref{rem:mu1}), for the $N$-Laplacian in $\R^{N}$, and for the biharmonic equation in $\R^{4}$ with critical growth in the sense of Adams \cite{RufSaniAdams13,SaniThesis}. Once the stationary picture is established, the dispersive consequences become natural directions of research. Examples are the orbital instability of the bound states and the role of a second $1$-node bound state in the dynamics above the ground-state threshold of \cite{DKM20}.

\section*{Declarations}
\noindent\textbf{Funding.} This research received no funding.

\noindent\textbf{Conflict of interest.} The author declares no conflict of interest.

\noindent\textbf{Data availability.} Data sharing is not applicable to this article, as no datasets were generated or analysed.

\end{document}